\documentclass{article}     
\usepackage{graphicx}
\usepackage{amssymb,amsmath}
\usepackage{setspace}

\newcommand{\qed}{\nobreak \ifvmode \relax \else
      \ifdim\lastskip<1.5em \hskip-\lastskip
      \hskip1.5em plus0em minus0.5em \fi \nobreak
      \vrule height0.75em width0.5em depth0.25em\fi}

\newtheorem{pro}{Proposition}
\newtheorem{thm}[pro]{Theorem}
\newtheorem{lem}[pro]{Lemma}
\newtheorem{cor}[pro]{Corollary}

\begin{document}

\title{Small time asymptotics for stochastic evolution equations}

\author{Terence Jegaraj \footnote{School of Mathematics and Statistics, UNSW, Sydney, NSW, 2052, Australia. Email: t.jegaraj@unsw.edu.au}}

\maketitle

\begin{abstract}
\noindent We obtain a large deviation principle describing the small time asymptotics of the solution of a stochastic evolution equation with multiplicative noise. Our assumptions are a condition on the linear drift operator that is satisfied by generators of analytic semigroups and Lipschitz continuity of the nonlinear coefficient functions. Methods originally used by Peszat~\cite{Pe} for the small noise asymptotics problem are adapted to solve the small time asymptotics problem. The results obtained in this way improve on some results of Zhang~\cite{Zh}.
\end{abstract}

\section{Introduction}\label{ps0}

\noindent In this paper we obtain a large deviation principle describing the small time asymptotics of the mild solution of a stochastic differential equation with Lipschitz continuous drift function $F$ and Lipschitz continuous diffusion function $G$, in a Hilbert space $H$:
\begin{equation}\label{pe4}
\left.\begin{array}{rcl}
dX(t)& = &(AX(t)+F(t,X(t)))\,dt + G(X(t))\,dW(t),\\
X(0)&=& x\in H.
\end{array}\right\}
\end{equation}
In equation (\ref{pe4}), $A$ is a linear operator generating a strongly continuous semigroup on $H$ and $(W(t))_{t\geq 0}$ is a cylindrical Wiener process. Peszat~\cite{Pe} found a large deviation principle to solve the corresponding small noise asymptotics problem, in which a small positive parameter $\epsilon$ multiplies the diffusion function $G$ and goes to zero. An important difference between this and the small time asymptotics problem is that, after changing variables in the small time asymptotics problem, a parameter $\epsilon$ multiplies the drift terms and so one must deal with stochastic convolutions of semigroups depending on $\epsilon$: $(S(\epsilon t):=\mathrm{e}^{\epsilon At})_{t\geq 0}$. Nevertheless Peszat's methods require little modification to yield a large deviation principle describing the behaviour as $\epsilon$ goes to zero of the continuous $H$-valued trajectories of $(X_x(\epsilon t))_{t\in[0,1]}$, where $(X_x(t))_{t\geq 0}$ is the solution of (\ref{pe4}).

\smallskip

\noindent In his pioneering paper~\cite{Va}, Varadhan studied the small time asymptotics of diffusion processes in $\mathbb{R}^n$. He obtained a large deviation principle describing the small time asymptotics of trajectories by starting from knowledge of the small time limiting behaviour of probability densities.

\smallskip

\noindent In the Hilbert space setting, Zhang~\cite{Zh} approached this problem using exponential equivalence arguments and obtained a large deviation principle describing the small time asymptotics of the mild solution of a stochastic equation with Lipschitz continuous and bounded diffusion function. To deal with the stochastic convolution term he assumed that the Hilbert space $H$ is compactly embedded in another Hilbert space $H_1$. Working in the space of continuous $H_1$-valued trajectories, Zhang was able to reduce the small time asymptotics problem to the case of zero drift, $A=0$ and $F=0$, where the small noise and small time asymptotics problems are equivalent. Zhang's result is a large deviation principle for distributions on the space of continuous $H_1$-valued trajectories, rather than the space of continuous $H$-valued trajectories. Using Peszat's methods, we avoid the need to introduce another Hilbert space corresponding to Zhang's $H_1$ and we also show that the large deviation principle holds if the diffusion function is not bounded.

\smallskip

\noindent In the next section we present some definitions, formulate our small time asymptotics problem precisely and list our main results. In succeeding sections we prove the main results.

\medskip

\section{Problem formulation and results}\label{ps1}

\noindent Let $(H,\langle\cdot,\cdot\rangle,|\cdot|)$ and $(U,\langle\cdot,\cdot\rangle_U,|\cdot|_U)$ be separable Hilbert spaces. Let $L(H,H)$ denote the space of bounded linear operators on $H$ and $L_2(U,H)$ denote the space of Hilbert-Schmidt operators from $U$ into $H$. Unless stated otherwise, we write $\|\cdot\|_E$ for the norm in any Banach space $E$. Let $A:D(A)\subset H\rightarrow H$ be the infinitesimal generator of a strongly continuous semigroup $(S(t))_{t\geq 0}$ of bounded linear operators  on $H$. Set
\[
M:=\sup_{t\in[0,1]}\|S(t)\|_{L(H,H)}.
\]
Let
\[
F:([0,1]\times H, \mathcal{B}_{[0,1]}\otimes\mathcal{B}_H)\rightarrow (H,\mathcal{B}_H)
\]
be a Borel measurable function and let  functions $F$ and
\[
G:H\rightarrow L_2(U,H)
\]
satisfy
\begin{eqnarray}
|F(t,x)-F(t,y)|&\le &\nu(t) |x-y|\,\,\,\mbox{ }\forall t\in[0,1] \mbox{ and }\forall x,y\in H \mbox{ and }\label{pe1a}\\
|F(t,x)|&\le &\nu(t)(1+|x|)\,\,\,\mbox{ }\forall t\in[0,1] \mbox{ and }\forall x\in H \mbox{ and }\label{pe1b}\\
\|G(x)-G(y)\|_{L_2(U,H)}&\le&\Lambda|x-y|\,\,\,\mbox{ }\forall x,y\in H \mbox{ and }\label{pe2a}\\
\|G(x)\|_{L_2(U,H)}&\le&\Lambda(1+|x|)\,\,\,\mbox{ }\forall x\in H,\label{pe2b}
\end{eqnarray}
where $\nu:[0,1]\rightarrow \mathbb{R}$ is a function in $L^2([0,1];\mathbb{R})$ and $\Lambda$ is a positive real constant.

\medskip

\noindent Let $(\Omega,\mathcal{F},P)$ be a probability space and let $(\mathcal{F}_t)_{t\geq 0}$ be a right continuous filtration of sub $\sigma$-algebras of $\mathcal{F}$ such that all sets in $\mathcal{F}$ of $P$ measure zero are in $\mathcal{F}_0$. Let $(g_k)$ be an orthonormal basis of $U$ and let $((\beta_k(t))_{t\geq 0})$ be an independent sequence of real valued $(\mathcal{F}_t)$-Brownian motions. A cylindrical Wiener process on $U$ is defined by the series
\[
W(t)=\sum_{k=1}^\infty \beta_k(t)\,g_k,
\]
which does not converge in $U$ but converges in an arbitrary Hilbert space $U_1$ containing $U$ and such that the embedding
\[
J:U\hookrightarrow U_1
\]
is Hilbert-Schmidt. Whatever our choice of $U_1$, the distribution of $W(1)$ in $U_1$ has reproducing kernel Hilbert space $U$. We now fix $U_1$ by taking a decreasing sequence of positive real numbers $(\lambda_k)$ such that $\sum_{k=1}^{\infty}\lambda_k^2<\infty$ and defining $U_1$ to be the completion of $U$ with the inner product
\begin{equation}\label{pe3}
\langle u,v\rangle_{U_1}:=\sum_{k=1}^\infty\lambda_k^2\langle u, g_k\rangle_U\langle v, g_k\rangle_U\,\,\,\mbox{ for all }u,v\in U.
\end{equation}
We abuse notation and denote the inner product on $U_1$ still by $\langle\cdot,\cdot\rangle_{U_1}$ and the norm on $U_1$ is denoted by $|\cdot|_{U_1}$.

\medskip

\noindent Our aim is to find a large deviation principle describing the small time asymptotics of the mild solution of the initial value problem in (\ref{pe4}).
The mild solution of (\ref{pe4}) is the $(\mathcal{F}_t)$-predictable process $(X_x(t))_{t\in[0,1]}$ such that
\begin{equation}\label{pe5a}
P\{\int_0^1|X_x(t)|^2\,dt<\infty\}=1
\end{equation}
and
\begin{equation}\label{pe5b}
X_x(t)=S(t)x+\int_0^t S(t-s)F(s,X_x(s))\,ds+\int_0^t S(t-s)G(X_x(s))\,dW(s)\,\,\,\mbox{ }P\mbox{ a.e.}
\end{equation}
for each $t\in [0,1]$.

\medskip

\noindent The existence, uniqueness and continuity result underlying this work is Theorem~\ref{ptA1} in the appendix of this paper.

\smallskip

\noindent Specifically, we will find a large deviation principle for the family of distributions on the space of continuous functions  mapping $[0,1]$ into $H$, $C([0,1];H)$:
\begin{equation}\label{pebefore9}
\mu^\epsilon_x:=\mathcal{L}(\,\omega\in\Omega\mapsto(t\in[0,1]\mapsto X_x(\epsilon t)(\omega))\,)\,\,\,\mbox{ }:\epsilon\in(0,1].
\end{equation}
From equation (\ref{pe5b}), for each $\epsilon\in(0,1]$ and $t\in[0,1]$ we have $P$ a.e.
\begin{eqnarray}
X_x(\epsilon t)
&=&S(\epsilon t)x+\int_0^{\epsilon t}S(\epsilon t-s)F(s,X_x(s))\,ds+\int_0^{\epsilon t}S(\epsilon t-s)G(X_x(s))\,dW(s)\nonumber\\
&=&S(\epsilon t)x+\epsilon\int_0^t S(\epsilon(t-u))F(\epsilon u,X_x(\epsilon u))\,du+\epsilon^{\frac{1}{2}}\int_0^t S(\epsilon(t-u))G(X_x(\epsilon u))\,dV^\epsilon(u),\mbox{\hspace{10mm}}\label{pe7}
\end{eqnarray}
where
\[
V^\epsilon(t):=\epsilon^{-\frac{1}{2}}W(\epsilon t)\,\,\,\mbox{ }\forall t\geq 0
\]
is a $U_1$-valued $(\mathcal{F}_{\epsilon t})$-Wiener process and $\mathcal{L}(V^\epsilon(1))=\mathcal{L}(W(1))$. By Proposition~\ref{ppA1}, for each $\epsilon\in(0,1]$ the continuous $(\mathcal{F}_t)$-predictable process $(X^\epsilon_x(t))_{t\in[0,1]}$  satisfying the equation
\begin{equation}\label{pe8}
X_x^\epsilon(t)=S(\epsilon t)x+\epsilon\int_0^tS(\epsilon(t-u))F(\epsilon u,X^\epsilon_x(u))\,du+\epsilon^\frac{1}{2}\int_0^tS(\epsilon(t-u))G(X^\epsilon_x(u))\,dW(u)
\end{equation}
$P$ a.e. for each $t\in[0,1]$ also has the distribution $\mu^\epsilon_x$ in trajectory space. Thus for each $\epsilon\in(0,1]$ we consider the process $(X_x^\epsilon(t))_{t\in[0,1]}$, which is the mild solution of the problem
\[
\left.
\begin{array}{rcl}
dX^\epsilon(t)&= &(\epsilon A X^\epsilon(t)+\epsilon F(\epsilon t,X^\epsilon(t)))\,dt+\epsilon^{\frac{1}{2}}G(X^\epsilon(t))\,dW(t)\\
X^\epsilon(0)&=&x
\end{array}\right\}
\]
and we define the corresponding trajectory-valued random variable
$X_x^\epsilon:\Omega\rightarrow C([0,1];H)$ by
\begin{equation}\label{pe9}
X^\epsilon_x(\omega):=(t\in[0,1]\mapsto X^\epsilon_x(t)(\omega))\,\,\,\mbox{ }\forall\omega\in\Omega.
\end{equation}

\bigskip

\noindent Recall that $\{\mu^\epsilon_x :\,\epsilon\in(0,1]\}$ is said to satisfy a large deviation principle with good rate function $\mathcal{I}_x: C([0,1];H)\rightarrow [0,\infty]$ if $\mathcal{I}_x$ has compact level set $\{u\in C([0,1];H):\,\mathcal{I}_x(u)\le r\}$ for each positive real $r$ and
\begin{eqnarray*}
\mbox{(i)}\;\;\;\liminf_{\epsilon\rightarrow 0}\epsilon\ln \mu^\epsilon_x(O)&\geq& -\inf_{u\in O}\mathcal{I}_x(u)\;\;\;\mbox{ for open subsets }O\mbox{ of } C([0,1];H)\mbox{ and }\\
\mbox{(ii)}\;\;\;\limsup_{\epsilon\rightarrow 0}\epsilon\ln\mu^\epsilon_x(C)&\le&-\inf_{u\in C}\mathcal{I}_x(u)\;\;\;\mbox{ for  closed subsets }C\mbox{ of } C([0,1];H).
\end{eqnarray*}
\noindent  We will use Freidlin's and Wentzell's (see~\cite[Theorem 3.3 on page 85]{FW} or~\cite[Proposition 12.2]{DPZ}) equivalent formulations of the lower bound condition, (i), and the upper bound condition, (ii). Before stating these, we introduce some more notation.

\noindent If $(E,\|\cdot\|_E)$ is a Banach space and $x$ is a point in $E$ then we denote the open ball in $E$ centred at $x$ and of radius $r>0$ by
$B_E(x,r):=\{e\in E:\|e-x\|_E<r\}$
and if $D$ is a subset of $E$ then we define $B_E(D,r):=\cup_{x\in D}B_E(x,r)$.

\noindent Freidlin and Wentzell showed that (i) and (ii) are, respectively, equivalent to: 
\begin{eqnarray*}
&&\!\!\!\!\!\!\!\!\!\!\mbox{(\textit{i})}\;\;\;\mbox{ for each }u\in C([0,1];H)\mbox{ and }\delta>0\mbox{ and }\gamma>0, \mbox{ there exists }\epsilon_0>0\mbox{ such that }\\
&&\mu^\epsilon_x(B_{C([0,1];H)}(u,\delta)) \geq \exp\left(\frac{-\mathcal{I}_x(u)-\gamma}{\epsilon}\right)\mbox{ for all }0<\epsilon<\epsilon_0\;\;\mbox{ and}\\
&&\!\!\!\!\!\!\!\!\!\!\mbox{(\textit{ii})}\;\;\;\mbox{ for each }r>0\mbox{ and }\delta>0\mbox{ and }\gamma>0, \mbox{ there exists }\epsilon_0>0\mbox{ such that }\\
&&\mu^\epsilon_x(B_{C([0,1];H)}(\{\mathcal{I}_x\le r\},\delta)) \geq  1-\exp\left(\frac{-r+\gamma}{\epsilon}\right)\mbox{ for all }0<\epsilon<\epsilon_0.
\end{eqnarray*}

\bigskip  

\noindent For each square integrable $U$-valued function $\phi\in L^2([0,1];U)$ and $x\in H$ we denote by $z^\phi_x$ the function in $C([0,1];H)$
such that
\[
z^\phi_x(t)=x+\int_0^tG(z^\phi_x(s))\phi(s)\,ds\,\,\,\mbox{ }\forall t\in[0,1].
\]
For each $x\in H$ we define the prospective rate function $\mathcal{I}_x:C([0,1];H)\rightarrow[0,\infty]$ by
\begin{equation}\label{pe10}
\mathcal{I}_x(u):=\frac{1}{2}\inf\left\{\int_0^1|\psi(s)|_U^2\,ds:\psi\in L^2([0,1];U)\mbox{ and }u=z^\psi_x\right\}
\end{equation}
for all $u\in C([0,1];H)$; here we take $\inf\emptyset =\infty$. We will prove the following theorem in Section~\ref{ps2}; it verifies that for each $x\in H$ the function $\mathcal{I}_x$ is well defined and a good rate function.
\begin{thm}\label{pt1}
\begin{enumerate}
\item Given $\phi\in L^2([0,1];U)$ and $x\in H$, $z^\phi_x$ is well defined; that is, there is a unique function $u\in C([0,1];H)$ such that
\[
u(t)=x+\int_0^tG(u(s))\phi(s)\,ds\,\,\,\mbox{ }\forall t\in[0,1].
\]

\item For fixed $u\in C([0,1];H)$ the linear operator
\[
\psi\in L^2([0,1];U)\mapsto \left(t\mapsto\int_0^tG(u(s))\psi(s)\,ds\right)\in C([0,1];H)
\]
is compact.
\item
Let $B\subset L^2([0,1];U)$ be weakly sequentially compact and let $K\subset H$ be compact. Then the set
\[
\mathcal{C}:=\{u\in C([0,1];H): u=z^\phi_x\mbox{ for some }\phi\in B\mbox{ and some }x\in K\}
\]
is compact. In particular $\{\mathcal{I}_x\le r\}$ is compact for any $x\in H$ and any $r\in (0,\infty)$ because the closed ball $\{\phi\in L^2([0,1];U):\|\phi\|_{L^2([0,1];U)}\le\sqrt{2r}\}$ is weakly sequentially compact.
\end{enumerate}
\end{thm}

\medskip

\noindent For each natural number $n$ let $\Pi_n:U\rightarrow U$ be the orthogonal projection of $U$ onto the span of $\{g_1,\ldots,g_n\}$:
\[
\Pi_nx:=\sum_{j=1}^n\langle x,g_j\rangle_U\, g_j\,\,\,\mbox{ }\forall x\in U.
\]
In our proof of the upper bound of the large deviation principle we use the fact that $\Pi_n$ can be written in terms of the bounded linear operator from $U_1$ into $U$:
\[
\Pi_n^1u:=\sum_{k=1}^n \lambda_k^{-2}\langle u,Jg_k\rangle_{U_1}\,g_k\,\,\,\mbox{ }\forall u\in U_1.
\]
We have $\Pi_n^1Jx=\Pi_nx$ for all $x\in U$, which follows from the definition of $U_1$.


\noindent The identity operator on any Banach space $E$ is denoted by $I_E$.
We can now state two additional assumptions (A1) and (A2) on $G$ and $(S(t))_{t\geq 0}$, respectively, which will only be used in the proof of the upper bound of the large deviation principle.

\noindent (\textbf{A1}) For each $r\in(0,\infty)$
\[
\sup_{h\in B_H(0,r)}\|G(h)(I_U-\Pi_n)\|_{L_2(U,H)}\rightarrow 0\,\,\,\mbox{ as }n\rightarrow\infty.
\]
(\textbf{A2}) For each $a\in(0,1]$ the family of functions in $L(H,H)$ with the norm topology:
\[
\{t\in[a,1]\mapsto S(\epsilon t)\in L(H,H)\,\,\,:\,\epsilon\in(0,1]\}
\]
\hspace{0mm}is uniformly equicontinuous.

\noindent Assumption (A1) is true when $G$ is of the form $G(x)=G_1(x)B\,\,\,\forall x\in H$, where $B$ is a constant operator in $L_2(U,U)$ and $G_1:H\rightarrow L(U,H)$ is bounded on bounded subsets of $H$.

\noindent Assumption (A2) is true when $(S(t))_{t\geq 0}$ is an analytic semigroup. Then there is a positive real constant $c$ such that $\|AS(t)\|_{L(H,H)}\le\frac{c}{t}$ for all $t\in(0,1]$ (see~\cite[Theorem 5.2 in chapter 2]{Pa}) and consequently $\|S(t)-S(r)\|_{L(H,H)}\le c\ln(\frac{t}{r})$ for all $t,r\in(0,1]$. We remark that in the small noise asymptotics paper~\cite{Pe} Peszat only needed to assume that $(S(t))_{t\geq 0}$ is continuous on $(0,\infty)$ in the norm topology.

\medskip
\noindent Our two main theorems are the following.
\begin{thm}\label{pt2}
Let $K$ be a compact subset of $H$ and let $\phi\in L^2([0,1];U)$. Let $\delta>0$ and $\gamma>0$. There exists $\epsilon_0>0$ such that for all $x\in K$ and for all $\epsilon\in(0,\epsilon_0]$
\[
P\left\{\sup_{t\in[0,1]}|X^\epsilon_x(t)-z^\phi_x(t)|<\delta\right\}\geq \exp\left(\frac{-\frac{1}{2}\int_0^1|\phi(s)|_U^2\,ds-\gamma}{\epsilon}\right).
\]
\end{thm}
\begin{thm}\label{pt3}
Assume that (A1) and (A2) hold. Let $K$ be a compact subset of $H$. Let $r>0$ and $\delta>0$ and $\gamma>0$. There exists $\epsilon_0>0$ such that for all $x\in K$ and for all $\epsilon\in(0,\epsilon_0]$
\[
P\{X^\epsilon_x\notin B_{C([0,1];H)}(\{\mathcal{I}_x\le r\},\,\delta)\}\le \exp\left(\frac{-r+\gamma}{\epsilon}\right).
\]
\end{thm}
The following result follows immediately from these theorems.
\begin{cor}\label{pc1}
Assume that (A1) and (A2) hold. Let $x\in H$. The family of distributions $\{\mu_x^\epsilon\,:\,\epsilon\in(0,1]\}$ defined in equation (\ref{pebefore9}) satisfies a large deviation principle with rate function $\mathcal{I}_x$.
\end{cor}
\textbf{Proof.} When $K=\{x\}$ Theorem~\ref{pt2} implies the Freidlin-Wentzell formulation of the lower bound of the large deviation principle of $\{\mathcal{L}(X^\epsilon_x)=\mu^\epsilon_x\,:\,\epsilon\in(0,1]\}$ with rate function $\mathcal{I}_x$ and Theorem~\ref{pt3} is the corresponding upper bound.$\qed$
\smallskip

\noindent
We will show in Section~\ref{ps3} that if Theorems~\ref{pt2} and~\ref{pt3} hold for bounded diffusion functions $G:H\rightarrow L_2(U,H)$ then the theorems also hold when  the function $G$ is not bounded.
Section~\ref{ps4} presents some important inequalities from Peszat's paper~\cite{Pe}, which are used to prove Theorems~\ref{pt2} and~\ref{pt3} in the case of bounded $G$ in Sections~\ref{ps5} and~\ref{ps6}.

\section{The rate function}\label{ps2}

\noindent In this section we prove Theorem~\ref{pt1}. Since $G$ is a Lipschitz continuous function, the proof of Theorem~\ref{pt1}(1) is a straightforward application of Banach's contraction mapping theorem and is omitted.

\medskip



\noindent\textbf{Proof of Theorem \ref{pt1}(2).} This proof follows the lines of the proof of~\cite[Proposition 8.4]{DPZ}. Let $u\in C([0,1];H)$. We want to show that the map
\[
\psi\in L^2([0,1];U)\mapsto \left(t\mapsto\int_0^tG(u(s))\psi(s)\,ds\right)\in C([0,1];H)
\]
is a compact linear operator.  We will show that an arbitrary bounded sequence $(\psi_n)$ in $L^2([0,1];U)$ is mapped to a sequence in $C([0,1];H)$ with a convergent subsequence.

\noindent Set $r:=\sup_{n\in\mathbb{N}}\|\psi_n\|_{L^2([0,1];U)}<\infty$. For $n\in\mathbb{N}$ and $0\le t<s\le 1$ we have
\begin{eqnarray*}
\left|\int_t^sG(u(\sigma))\psi_n(\sigma)\,d\sigma\right|
&\le& r\sup_{\sigma\in[0,1]}\|G(u(\sigma))\|_{L_2(U,H)}\sqrt{s-t}.
\end{eqnarray*}
Thus the family of functions
\[
\left\{t\in[0,1]\mapsto\int_0^tG(u(s))\psi_n(s)\,ds\in H\,\,\,:\,\,n\in \mathbb{N}\right\}
\]
is uniformly equicontinuous. We will show that there is a subsequence of the sequence of continuous functions
\[
\left(t\in[0,1]\mapsto\int_0^tG(u(s))\psi_n(s)\,ds\in H\right)
\]
which converges pointwise on a dense subset of $[0,1]$; then, as this subsequence is uniformly equicontinuous, it is Cauchy in $C([0,1];H)$ and we will be done.

\noindent For each $t\in(0,1]$ define the linear operator $A_t:L^2([0,1];U)\rightarrow H$ by
\[
A_t\psi=\int_0^tG(u(s))\psi(s)\,ds\,\,,\,\,\,\mbox{ }\psi\in L^2([0,1];U).
\]
One can show by direct computation that for each $t\in(0,1]$ $A_t$ is Hilbert-Schmidt. For each $t_i\in(0,1]\cap\mathbb{Q}$, since $A_{t_i}$ is a compact operator, the set $\{A_{t_i}\psi_n\,\,:\,\,\,n\in\mathbb{N}\}$ is relatively compact in $H$. We can apply the diagonal argument to the sequence of sequences
\[
\begin{array}{ccccc}
A_{t_1}\psi_1& A_{t_1}\psi_2& A_{t_1}\psi_3& A_{t_1}\psi_4&\cdots\\
A_{t_2}\psi_1& A_{t_2}\psi_2& A_{t_2}\psi_3& A_{t_2}\psi_4&\cdots\\
\vdots&\vdots&\vdots&\vdots&{}\\
A_{t_i}\psi_1& A_{t_i}\psi_2& A_{t_i}\psi_3& A_{t_i}\psi_4&\cdots\\
\vdots&\vdots&\vdots&\vdots&{}
\end{array}
\]
to conclude that there is a strictly increasing sequence of natural numbers $(n_k)$ such that
\[
\lim_{k\rightarrow\infty}A_{t_i}\psi_{n_k}\,\,\,\mbox{ exists for each }i\in\mathbb{N}.\,\,\,\,\,\,\,\,\qed
\]

\medskip

\noindent\textbf{Proof of Theorem \ref{pt1}(3).} Let $B\subset L^2([0,1];U)$ be weakly sequentially compact and let $K\subset H$ be compact. We want to show that
\[
\mathcal{C}:=\{u\in C([0,1];H)\,:\,\,u=z^\phi_x\,\,\,\mbox{ for some }\,\,\phi\in B\,\,\mbox{ and some }\,\,x\in K\}
\]
is compact. Weak sequential compactness of $B$ implies that $B$ is bounded. Hence
\[
q:=\sup\{\|\psi\|_{L^2([0,1];U)}\,:\,\psi\in B\}
\]
is finite. Let $(u_n)$ be a sequence of elements of $\mathcal{C}$. For each $n\in\mathbb{N}$ there is $\phi_n\in B$ and $x_n\in K$ such that $u_n=z^{\phi_n}_{x_n}$. By compactness of $K$ and weak sequential compactness of $B$, there is a strictly increasing sequence of natural numbers $(n_k)$ and there are vectors $x\in K$ and $\phi\in B$ such that $x_{n_k}$ converges to $x$ in $H$ and $\phi_{n_k}$ converges to $\phi$ in the weak topology of $L^2([0,1];U)$ as $k$ goes to infinity. We claim that $u_{n_k}\rightarrow u:=z^\phi_x$ as $k\rightarrow\infty$. Indeed by Gronwall's Lemma we have
\[
\sup_{t\in[0,1]}|u(t)-u_{n_k}(t)|\le\left(|x-x_{n_k}|+\sup_{r\in[0,1]}\left|\int_0^rG(u(s))(\phi(s)-\phi_{n_k}(s))\,ds\right|\right)\exp(\Lambda q) 
\]
for each $k\in\mathbb{N}$ and the right hand side of the above inequality goes to zero as $k\rightarrow\infty$ because of compactness of the linear operator
\[
\psi\in L^2([0,1];U)\mapsto\left(t\mapsto\int_0^tG(u(s))\psi(s)\,ds\right)\in C([0,1];H).\,\,\,\,\,\,\qed
\]

\section{Reducing the problem to the case of bounded $G$}\label{ps3}

\noindent In this section we show that if Theorems~\ref{pt2} and~\ref{pt3} hold under the additional assumption that the function $G:H\rightarrow L_2(U,H)$ is bounded then they hold also for $G$ which is not bounded. This idea is copied from Cerrai and R$\ddot{\mathrm{o}}$ckner~\cite[Theorem 6.4]{CR}.

\medskip

\noindent For each $R\in(0,\infty)$ define $G_R:H\rightarrow L_2(U,H)$ by
\[
G_R(x):=\left\{\begin{array}{ll}
G(x)&\mbox{ if }|x|\le R\\
G(\textstyle{\frac{R}{|x|}}x)&\mbox{ if }|x|>R;
\end{array}\right.
\]
it is straightforward to show that $\sup_{x\in H}\|G_R(x)\|_{L_2(U,H)}<\infty$ and that inequalities~(\ref{pe2a}) and~(\ref{pe2b}) also hold with $G_R$ in place of $G$.
For each $x\in H$ and $R\in(0,\infty)$ define $\mathcal{I}_{R,x}:C([0,1];H)\rightarrow[0,\infty]$ by
\[\begin{array}{l}
\mathcal{I}_{R,x}(u):=\\
\,\,\,\,\,\frac{1}{2}\inf\left\{\int_0^1|\phi(s)|^2_U\,ds\,:\,\phi\in L^2([0,1];U) \mbox{ and }u(t)=x+\int_0^tG_R(u(s))\phi(s)\,ds\,\,\,\mbox{ }\forall t\in[0,1]\right\},
\end{array}\]
for all $u\in C([0,1];H)$; here we take $\inf\emptyset=\infty$.

\noindent For each $x\in H$, $R\in(0,\infty)$ and $\epsilon\in(0,1]$ define $(X^\epsilon_{R,x}(t):(\Omega,\mathcal{F}_t)\rightarrow(H,\mathcal{B}_H))_{t\in[0,1]}$ to be the continuous $(\mathcal{F}_t)$-predictable process satisfying
\begin{equation}\label{pe11}
X^\epsilon_{R,x}(t)=S(\epsilon t)x+\epsilon\int_0^tS(\epsilon(t-s))F(\epsilon s,X^\epsilon_{R,x}(s))\,ds+\epsilon^\frac{1}{2}\int_0^tS(\epsilon(t-s))G_R(X^\epsilon_{R,x}(s))\,dW(s)
\end{equation}
$P$ a.e. for each $t\in[0,1]$ and let $X^\epsilon_{R,x}:\Omega\rightarrow C([0,1];H)$ be the corresponding trajectory-valued random variable:
\[
X^\epsilon_{R,x}(\omega):=(t\mapsto X^\epsilon_{R,x}(t)(\omega))\,\,\,\,\,\mbox{ }\forall \omega\in\Omega.
\]
Recall that we defined $(X^\epsilon_x(t))_{t\in[0,1]}$ and $X^\epsilon_x$ in equations~(\ref{pe8}) and~(\ref{pe9}).

\begin{lem}\label{pl4}
Let $\rho\in(0,\infty)$. Given $r\in(0,\infty)$ and $\delta\in(0,\infty)$ there exists $R\in(0,\infty)$ such that
\begin{enumerate}
\item for each $x\in B_H(0,\rho)$
\[
\{\mathcal{I}_{R,x}\le r\}=\{\mathcal{I}_x\le r\}
\]
and
\item for each $x\in B_H(0,\rho)$ and for each $\epsilon\in(0,1]$
\[
P\{X^\epsilon_x\in B_{C([0,1];H)}(\{\mathcal{I}_x\le r\},\delta)\}=P\{X^\epsilon_{R,x}\in B_{C([0,1];H)}(\{\mathcal{I}_x\le r\},\delta)\}.
\]
\end{enumerate}
\end{lem}
\textbf{Proof.} Let $r>0$ and $\delta>0$. Set
\begin{equation}\label{pe12}
R:=(\rho+\Lambda\sqrt{2r})\exp(\Lambda\sqrt{2r})+\delta.
\end{equation}
Let $x\in B_H(0,\rho)$.

\smallskip

\noindent We firstly prove part (1).

\noindent Suppose $u\in C([0,1];H)$ and $\mathcal{I}_x(u)\le r$. Then there exists $\phi\in L^2([0,1];U)$ such that $\|\phi\|_{L^2([0,1];U)}\le\sqrt{2r}$ and
\begin{equation}\label{pe13a}
u(t)=x+\int_0^tG(u(s))\,\phi(s)\,ds\,\,\,\mbox{ }\forall t\in[0,1].
\end{equation}
Taking the norm of both sides of this equation and applying Gronwall's Lemma gives
\begin{equation}\label{pe13}
\sup_{t\in[0,1]}|u(t)|<(\rho+\Lambda\sqrt{2r})\exp(\Lambda\sqrt{2r}).
\end{equation}
Since $G(x)=G_R(x)$ for all $x\in B_H(0,R)$, (\ref{pe13}) implies that $u$ satisfies
\[
u(t)=x+\int_0^tG_R(u(s))\phi(s)\,ds\,\,\,\mbox{ }\forall t\in[0,1]
\]
and $\mathcal{I}_{R,x}(u)\le\frac{1}{2}\int_0^1|\phi(s)|_U^2\,ds\le r$. We have shown that $\{\mathcal{I}_x\le r\}\subset\{\mathcal{I}_{R,x}\le r\}$ and the reverse inclusion is proved in a similar way.

\medskip

\noindent We now prove part (2). Let $\epsilon\in(0,1]$.

\noindent Define the $(\mathcal{F}_t)$-stopping time
\[
\tau(\omega):=\inf\{t\in[0,1]:|X^\epsilon_x(t)(\omega)|\geq R\}\,\,,\,\,\,\mbox{ }\omega\in\Omega,
\]
where we take $\tau(\omega)=1$ if $|X^\epsilon_x(t)(\omega)|<R$ for all $t\in[0,1]$. By our choice of $R$ and inequality (\ref{pe13}) we have
\[
B_{C([0,1];H)}(\{\mathcal{I}_x\le r\},\,\delta)\subset B_{C([0,1];H)}(0,R).
\]
Thus we can tell if the trajectory $X^\epsilon_x(\omega)$ lies in $B_{C([0,1];H)}(\{\mathcal{I}_x\le r\},\,\delta)$ by observing the trajectory just up to time $\tau(\omega)$; also for $P$ a.e. $\omega\in\Omega$ we have $\sup_{t\in[0,\tau(\omega)]}|X^\epsilon_x(t)(\omega)|\le R$ and if $\tau(\omega)<1$ then $|X^\epsilon_x(\tau(\omega))(\omega)|=R$.

\noindent Let $t\in(0,1]$. We have
\[
X^\epsilon_x(t)=S(\epsilon t)x+\epsilon\int_0^tS(\epsilon(t-s))F(\epsilon s,X^\epsilon_x(s))\,ds+\epsilon^\frac{1}{2}\int_0^tS(\epsilon(t-s))G(X^\epsilon_x(s))\,dW(s)\,\,\,\,\,\mbox{ }P\mbox{ a.e..}
\]
Multiplying both sides of this equation by the indicator of the stochastic interval $[0,\tau]$ we have
\begin{eqnarray*}
1_{[0,\tau]}(t)X^\epsilon_x(t)&=& 1_{[0,\tau]}(t)S(\epsilon t)x+1_{[0,\tau]}(t)\epsilon\int_0^tS(\epsilon(t-s))F(\epsilon s,1_{[0,\tau]}(s)X^\epsilon_x(s))\,ds+\\
&&\mbox{\hspace{3mm}}1_{[0,\tau]}(t)\epsilon^\frac{1}{2}\int_0^{t\wedge\tau}1_{[0,t]}(s)S(\epsilon(t-s))G(X^\epsilon_x(s))\,dW(s)\,\,\,\,\,\mbox{ }P\mbox{ a.e.}\\
&=&1_{[0,\tau]}(t)S(\epsilon t)x+1_{[0,\tau]}(t)\epsilon\int_0^tS(\epsilon(t-s))F(\epsilon s,1_{[0,\tau]}(s)X^\epsilon_x(s))\,ds+\\
&&\mbox{\hspace{3mm}}1_{[0,\tau]}(t)\epsilon^\frac{1}{2}\int_0^t1_{[0,\tau]}(s)1_{[0,t]}(s)S(\epsilon(t-s))G_R(1_{[0,\tau]}(s)X^\epsilon_x(s))\,dW(s)
\end{eqnarray*}
$P$ a.e..

\noindent From~(\ref{pe11}) we can obtain a similar equality for $1_{[0,\tau]}(t)X^\epsilon_{R,x}(t)$. Therefore
\begin{eqnarray}
&&1_{[0,\tau]}(t)(X^\epsilon_x(t)-X^\epsilon_{R,x}(t))=\mbox{\hspace{10cm}}\nonumber\\
&&\mbox{\hspace{5mm}}1_{[0,\tau]}(t)\epsilon\int_0^tS(\epsilon(t-s))[F(\epsilon s,1_{[0,\tau]}(s)X^\epsilon_x(s))-F(\epsilon s,1_{[0,\tau]}(s)X^\epsilon_{R,x}(s))]ds+\mbox{\hspace{5mm}}\nonumber\\
&&\mbox{\hspace{5mm}}1_{[0,\tau]}(t)\epsilon^\frac{1}{2}\int_0^t1_{[0,\tau]}(s)S(\epsilon(t-s))[G_R(1_{[0,\tau]}(s)X^\epsilon_x(s))-G_R(1_{[0,\tau]}(s)X^\epsilon_{R,x}(s))]\,dW(s)\mbox{\hspace{0mm}}\label{pe14}
\end{eqnarray}
$P$ a.e..

\noindent By Theorem~\ref{ptA1} we have
\[
\sup_{u\in[0,1]}E\left[|X^\epsilon_x(u)|^2\right]<\infty\,\,\,\mbox{ and }\,\,\,\sup_{u\in[0,1]}E\left[|X^\epsilon_{R,x}(u)|^2\right]<\infty.
\]
Thus, taking norms on both sides of equation (\ref{pe14}), then squaring both sides and taking expectations, we obtain
\begin{eqnarray}
&&E\left[|1_{[0,\tau]}(t)(X^\epsilon_x(t)-X^\epsilon_{R,x}(t))|^2\right]\mbox{\hspace{9cm}}\nonumber\\
&&\le\,\,2\epsilon M^2(\textstyle{\int_0^1\nu^2(s)\,ds}+\Lambda^2)\int_0^tE\left[|1_{[0,\tau]}(s)(X^\epsilon_x(s)-X^\epsilon_{R,x}(s))|^2\right]\,ds\,\,\,\,\,\mbox{ for each } t\in[0,1].\nonumber
\end{eqnarray}
By Gronwall's Lemma it follows that
\[
1_{[0,\tau]}(t)X^\epsilon_x(t)=1_{[0,\tau]}(t)X^\epsilon_{R,x}(t)\,\,\,\mbox{ }\forall t\in[0,1]\,\,\,\mbox{ }P\mbox{ a.e..}
\]
We conclude that for $P$ a.e. $\omega\in\Omega$,
\begin{enumerate}
\item if $\tau(\omega)=1$ then $X^\epsilon_x(\omega)=X^\epsilon_{R,x}(\omega)$ and
\item if $\tau(\omega)<1$ then $|X^\epsilon_x(\tau(\omega))(\omega)|=|X^\epsilon_{R,x}(\tau(\omega))(\omega)|=R$ and trajectories $X^\epsilon_x(\omega)$ and $X^\epsilon_{R,x}(\omega)$ do not belong to $B_{C([0,1];H)}(\{\mathcal{I}_x\le r\},\,\delta)$.
\end{enumerate}
Thus
\[
P\{X^\epsilon_x\in B_{C([0,1];H)}(\{\mathcal{I}_x\le r\},\delta)\}=P\{X^\epsilon_{R,x}\in B_{C([0,1];H)}(\{\mathcal{I}_x\le r\},\delta)\}.\,\,\,\,\,\,\,\qed
\]

\bigskip

\noindent Given $x\in H$ and $\phi\in L^2([0,1];U)$ and $R>0$ we denote by $z^\phi_{R,x}$ the function $u\in C([0,1];H)$ such that $u(t)=x+\int_0^tG_R(u(s))\phi(s)\,ds$ for all $t\in[0,1]$; recall that $z^\phi_x$ is the function $v\in C([0,1];H)$ such that $v(t)=x+\int_0^tG(v(s))\phi(s)\,ds$ for all $t\in[0,1]$. The proof of the next lemma is similar to that of Lemma~\ref{pl4} and is omitted.
\begin{lem}\label{pl5}
Let $K\subset H$ be compact. Given $\phi\in L^2([0,1];U)$ and $\delta>0$ there exists $R>0$ such that:
\begin{enumerate}
\item for all $x\in K$ we have
\[
z^\phi_x=z^\phi_{R,x}
\]
and
\item for all $x\in K$ and all $\epsilon\in(0,1]$ we have
\[
P\{X^\epsilon_x\in B_{C([0,1];H)}(z^\phi_x,\,\delta)\}=P\{X^\epsilon_{R,x}\in B_{C([0,1];H)}(z^\phi_x,\,\delta)\}.
\]
\end{enumerate}
\end{lem}
Lemmas~\ref{pl4} and~\ref{pl5} have the following important corollaries.
\begin{cor}\label{pc6}
Suppose that Theorem~\ref{pt2} holds under the additional assumption that the diffusion function $G:H\rightarrow L_2(U,H)$ is bounded. Then it also holds if the function $G$ is not bounded.
\end{cor}
\begin{cor}\label{pc7}
Suppose that Theorem~\ref{pt3} holds under the additional assumption that the diffusion function $G:H\rightarrow L_2(U,H)$ is bounded. Then it also holds if the function $G$ is not bounded.
\end{cor}
Thanks to Corollaries~\ref{pc6} and~\ref{pc7} our task reduces to proving Theorems~\ref{pt2} and~\ref{pt3} under the additional assumption:

\medskip

\noindent(\textbf{A3}) the function $G:H\rightarrow L_2(U,H)$ is bounded, that is:
\[
\sup_{x\in H}\|G(x)\|_{L_2(U,H)}<\infty.
\]

\section{Exponential bounds}\label{ps4}

\noindent To prove Theorems~\ref{pt2} and~\ref{pt3} in the case of bounded $G$ we shall need some exponential tail estimates for stochastic integrals and stochastic convolutions due to Chow and Menaldi~\cite{CM} and Peszat~\cite{Pe}. The formulations we present without proof are Peszat's~\cite{Pe}.

\smallskip

\noindent Let $\mathcal{P}_1$ denote the $(\mathcal{F}_t)$-predictable $\sigma$-algebra of $[0,1]\times\Omega$. Let $\xi:([0,1]\times\Omega,\mathcal{P}_1)\rightarrow (L_2(U,H),\mathcal{B}_{L_2(U,H)})$
be a measurable function.

\begin{thm}[Chow's and Menaldi's bound for stochastic integrals]\label{pt8}
If there exists a positive real number $\eta_1$ such that
\[
\int_0^1\|\xi(s)\|_{L_2(U,H)}^2\,ds\le \eta_1\,\,\,\,\,\mbox{ }P\mbox{ a.e.}
\]
then for any $\delta>0$
\[
P\left\{\sup_{t\in[0,1]}\left|\int_0^t\xi(s)\,dW(s)\right|\geq\delta\right\}\le\,3\exp\left(-\frac{\delta^2}{4\eta_1}\right).
\]
\end{thm}

\begin{thm}[Peszat's bound for stochastic convolutions]\label{pt9}
Let $(T(t))$ be a strongly continuous semigroup of bounded linear operators on $H$. Suppose $\alpha_0\in(0,\frac{1}{2})$ and $p_0>1$ are such that
\[
\kappa:=\left(\int_0^1t^{(\alpha_0-1)p_0}\|T(t)\|_{L(H,H)}^{p_0}\,dt\right)^\frac{1}{p_0}<\infty.
\]
If there exists a positive real number $\eta_2$ such that
\[
\sup_{t\in[0,1]}\int_0^t(t-s)^{-2\alpha_0}\|T(t-s)\xi(s)\|_{L_2(U,H)}^2\,ds\le\eta_2\,\,\,\,\,\mbox{ }P\mbox{ a.e.}
\]
then the process $(\int_0^tT(t-s)\xi(s)\,dW(s))_{t\in[0,1]}$ has a continuous version in $H$ and for any $\delta>0$
\[
P\left\{\sup_{t\in[0,1]}\left|\int_0^tT(t-s)\xi(s)\,dW(s)\right|\geq\delta\right\}\le\,C\exp\left(-\frac{\delta^2}{\kappa^2\eta_2}\right)
\]
where $C=4+\exp(4n_0!)^\frac{1}{n_0}$ and $n_0=\frac{p_0}{2p_0-2}+1$.
\end{thm}
In the proof of Theorem~\ref{pt3} we also use a large deviation principle associated with the trajectory-valued random variable $W:(\Omega,\mathcal{F},P)\rightarrow(C([0,1];U_1),\mathcal{B}_{C([0,1];U_1)})$ defined by
\[
W(\omega):=(t\in[0,1]\mapsto W(t)(\omega)\in U_1)\,\,\,\,\,\mbox{ }\forall\omega\in\Omega.
\]
As shown in~\cite[Theorem 1 in Section 6.2]{Za}, the distribution of $W$ is symmetric Gaussian and its reproducing kernel Hilbert space is
\[
H_W:=\left\{t\in[0,1]\mapsto J\int_0^t\psi(s)\,ds\,:\,\psi\in L^2([0,1];U)\right\},
\]
with norm $\|\cdot\|_{H_W}$ defined by
\[
\|y\|_{H_W}^2:=\int_0^1|\psi(s)|_U^2\,ds\,:\,\,\,\psi\in L^2([0,1];U)\mbox{ and }y(t)=J\int_0^t\psi(s)\,ds\,\,\,\mbox{ }\forall t\in[0,1].
\]
Thus, by~\cite[Theorem 12.7]{DPZ}, the family of Gaussian measures
\[
\{\mathcal{L}(\epsilon^\frac{1}{2}W\,:\,(\Omega,\mathcal{F},P)\rightarrow(C([0,1];U_1),\mathcal{B}_{C([0,1];U_1)}))\,\,\,:\,\,\,\epsilon\in(0,1]\}
\]
satisfies a large deviation principle with rate function $\mathcal{I}_W:C([0,1];U_1)\rightarrow [0,\infty]$ defined by
\begin{equation}\label{pe19}
\mathcal{I}_W(f):=\left\{\begin{array}{cl}
\frac{1}{2}\|f\|^2_{H_W}& \mbox{ if }f\in H_W,\\
\infty&\mbox{ if }f\notin H_W.
\end{array}\right.
\end{equation}

\section{The lower bound}\label{ps5}

In this section we prove Theorem~\ref{pt2} under the additional assumption (A3).

\noindent\textbf{Proof of Theorem~\ref{pt2} assuming (A3).} Let $K\subset H$ be compact and fix $\phi\in L^2([0,1];U)$. Recall that for each $x\in K$ $z^\phi_x\in C([0,1];H)$ satisfies
\[
z^\phi_x(t)=x+\int_0^tG(z^\phi_x(s))\phi(s)\,ds\,\,\,\mbox{ }\forall t\in[0,1].
\]
Fix $\delta>0$ and $\gamma>0$. For each $\epsilon\in(0,1]$ define the process $(W^\epsilon(t):(\Omega,\mathcal{F}_t)\rightarrow (U_1,\mathcal{B}_{U_1}))_{t\in[0,1]}$ by
\begin{equation}\label{pe20}
W^\epsilon(t):=W(t)-\epsilon^{-\frac{1}{2}}J\int_0^t\phi(s)\,ds\,\,\,\mbox{ }\forall t\in[0,1].
\end{equation}
By~\cite[Theorem 10.14]{DPZ} $(W^\epsilon(t))_{t\in[0,1]}$ is a Wiener process with respect to filtration $(\mathcal{F}_t)$ on probability space $(\Omega,\mathcal{F},P^\epsilon)$ where
\begin{equation}\label{peafter24}
dP^\epsilon(\omega)=\exp\left(\epsilon^{-\frac{1}{2}}\int_0^1\langle\phi(s),\cdot\rangle_U\,dW(s)(\omega)-\frac{1}{2\epsilon}\int_0^1|\phi(s)|_U^2\,ds\right)\,dP(\omega)
\end{equation}
and $P^\epsilon(W^\epsilon(1))^{-1}=P(W(1))^{-1}$.

\noindent Taking the reciprocal of the Radon-Nikodym derivative in equation~(\ref{peafter24}) and using Lemma~\ref{pAl1} to replace the It$\hat{\mathrm{o}}$ integral on the right hand side by one with respect to $(W^\epsilon(t))_{t\in[0,1]}$:
\[
\int_0^1\langle\phi(s),\cdot\rangle_U\,dW^\epsilon(s)=\int_0^1\langle\phi(s),\cdot\rangle_U\,dW(s)-\epsilon^{-\frac{1}{2}}\int_0^1|\phi(s)|^2_U\,ds\,\,\,\,\,\mbox{ }P^\epsilon\mbox{ a.e.,}
\]
we have
\[
dP(\omega)=\exp\left(-\epsilon^{-\frac{1}{2}}\int_0^1\langle\phi(s),\cdot\rangle_U\,dW^\epsilon(s)(\omega)-\frac{1}{2\epsilon}\int_0^1|\phi(s)|^2_U\,ds\right)\,dP^\epsilon(\omega).
\]
To shorten notation, for each $x\in K$ and each $\epsilon\in(0,1]$ set
\begin{eqnarray*}
\mathcal{A}(\epsilon,x)&:=&\left\{\omega\in\Omega:\sup_{t\in[0,1]}|X^\epsilon_x(t)(\omega)-z^\phi_x(t)|<\delta\right\}\,\,\,\,\,\mbox{ and }\\
\mathcal{D}(\epsilon)&:=&\left\{\omega\in\Omega:\,\left|\epsilon^\frac{1}{2}\int_0^1\langle\phi(t),\cdot\rangle_U\,dW^\epsilon(t)(\omega)\right|\le\frac{\gamma}{2}\right\}.
\end{eqnarray*}
We have
\begin{eqnarray*}
P(\mathcal{A}(\epsilon,x))&=&\int_{\Omega}1_{\mathcal{A}(\epsilon,x)}\,\exp\left(-\epsilon^{-\frac{1}{2}}\int_0^1\langle\phi(s),\cdot\rangle_U\,dW^\epsilon(s)-\frac{1}{2\epsilon}\int_0^1|\phi(s)|_U^2\,ds\right)\,dP^\epsilon\\
&\geq&\int_\Omega 1_{\mathcal{A}(\epsilon,x)\cap\mathcal{D}(\epsilon)}\exp\left(-\epsilon^{-\frac{1}{2}}\int_0^1\langle\phi(s),\cdot\rangle_U\,dW^\epsilon(s)-\frac{1}{2\epsilon}\int_0^1|\phi(s)|_U^2\,ds\right)\,dP^\epsilon\\
&\geq&\exp\left(-\frac{\gamma}{2\epsilon}-\frac{1}{2\epsilon}\int_0^1|\phi(s)|^2_U\,ds\right)P^\epsilon(\mathcal{A}(\epsilon,x)\cap\mathcal{D}(\epsilon)).
\end{eqnarray*}
It remains to show that there exists $\epsilon_0>0$ such that $P^\epsilon(\mathcal{A}(\epsilon,x)\cap\mathcal{D}(\epsilon))\geq\exp(-\frac{\gamma}{2\epsilon})$  for all $x\in K$ and for all $\epsilon\in(0,\epsilon_0]$. We will actually show something more:
$
P^\epsilon(\mathcal{A}(\epsilon,x)^c\cup\mathcal{D}(\epsilon)^c)\rightarrow 0\,\,\,\mbox{ as }\,\,\,\epsilon\rightarrow 0\,\,\,\mbox{ uniformly in }x\in K.
$

\smallskip

\noindent Let $\epsilon\in(0,1]$ and let $x\in H$.

\noindent For each $t\in(0,1]$ we have
\begin{eqnarray}
&&\!\!\!\!\!\!\!\!|X^\epsilon_x(t)-z^\phi_x(t)|\nonumber\\
&=&\left|S(\epsilon t)x-x+\epsilon\int_0^tS(\epsilon(t-s))(F(\epsilon s,X^\epsilon_x(s))-F(\epsilon s,z^\phi_x(s)))\,ds\right.\nonumber\\
&&{}+\,\,\epsilon\int_0^tS(\epsilon(t-s))F(\epsilon s,z^\phi_x(s))\,ds\nonumber\\
&&{}+\int_0^tS(\epsilon(t-s))(G(X^\epsilon_x(s))-G(z^\phi_x(s)))\phi(s)\,ds\nonumber\\
&&{}+\int_0^t(S(\epsilon(t-s))-I_H)G(z^\phi_x(s))\phi(s)\,ds\nonumber\\
&&\left.{}+\epsilon^{\frac{1}{2}}\left(\int_0^tS(\epsilon(t-s))G(X^\epsilon_x(s))\,dW(s)-\epsilon^{-\frac{1}{2}}\int_0^tS(\epsilon(t-s))G(X^\epsilon_x(s))\phi(s)\,ds\right)\right|\nonumber\\
&\le&\sup_{r\in[0,\epsilon]}|S(r)x-x|+\epsilon M\int_0^t\nu(\epsilon s)|X^\epsilon_x(s)-z^\phi_x(s)|\,ds\nonumber\\
&&{}+\epsilon M\int_0^1\nu(\epsilon s)(1+|z^\phi_x(s)|)\,ds\nonumber\\
&&{}+M\Lambda\int_0^t|X^\epsilon_x(s)-z^\phi_x(s)||\phi(s)|_U\,ds\nonumber\\
&&{}+\sup\{\|(S(r)-I_H)G(z^\phi_x(s))\|_{L_2(U,H)}\,:\,r\in[0,\epsilon]\,,\,\,\,s\in[0,1]\}\left(\int_0^1|\phi(s)|_U^2\,ds\right)^\frac{1}{2}\nonumber\\
&&{}+\sup_{r\in[0,1]}\epsilon^\frac{1}{2}\left|\int_0^rS(\epsilon(r-s))G(X^\epsilon_x(s))\,dW^\epsilon(s)\right|\,\,\,\,\,\mbox{ }P^\epsilon\mbox{ a.e.}.\mbox{\hspace{10mm}}\label{pe21}
\end{eqnarray}
The last term on the right of (\ref{pe21}) is obtained by applying Lemma~\ref{pAl1}.
Squaring both sides of inequality~(\ref{pe21}) and applying Gronwall's Lemma we have
\begin{eqnarray}
&&\!\!\!\!\!\!\!\!\sup_{t\in[0,1]}|X^\epsilon_x(t)-z^\phi_x(t)|^2\mbox{\hspace{8cm}}\nonumber\\
&\le&6\left[\sup_{r\in[0,\epsilon]}|S(r)x-x|^2+\epsilon M^2(\textstyle{\int_0^1\nu^2(s)\,ds})\int_0^1(1+|z_x^\phi(s)|)^2\,ds+{}\right.\nonumber\\
&&\,\,\,\,\,\sup\{\|(S(r)-I_H)G(z^\phi_x(s))\|_{L_2(U,H)}\,:\,r\in[0,\epsilon]\,,\,\,\,s\in[0,1]\}^2\int_0^1|\phi(s)|_U^2\,ds+{}\nonumber\\
&&\,\,\,\,\,\left.\sup_{r\in[0,1]}\left|\epsilon^\frac{1}{2}\int_0^rS(\epsilon(r-s))G(X^\epsilon_x(s))\,dW^\epsilon(s)\right|^2\right]\nonumber\\
&&\mbox{\hspace{20mm}}{}\times\exp\left(6M^2\left(\int_0^1\nu^2(s)\,ds+\Lambda^2\int_0^1|\phi(s)|_U^2\,ds\right)\right)\,\,\,\mbox{ }P^\epsilon\mbox{ a.e..}\,\,\,\,\,\,\,\,\,\,\,\,\,\,\,\,\,\,\,\label{e24a}
\end{eqnarray}
Since $K$ is compact, there exists $\epsilon_1>0$ such that for all $x\in K$ and for all $\epsilon\in(0,\epsilon_1]$ we have
\begin{eqnarray}
&&\!\!\!\!\!\!\!\!\!\!\!\!P^\epsilon\{\sup_{t\in[0,1]}|X^\epsilon_x(t)-z^\phi_x(t)|\geq\delta\}\mbox{\hspace{8cm}}\nonumber\\
&\le&P^\epsilon\left\{\sup_{r\in[0,1]}\left|\textstyle{\int_0^rS(\epsilon(r-s))G(X^\epsilon_x(s))\,dW^\epsilon(s)}\right|\geq\frac{\delta}{3\epsilon^\frac{1}{2}\exp(3M^2(\textstyle{\int_0^1\nu^2(s)\,ds}+\Lambda^2\textstyle{\int_0^1|\phi(s)|_U^2\,ds}))}\right\}.\mbox{\hspace{7mm}}\label{pe23}
\end{eqnarray}
Since we are assuming (A3), we can apply Peszat's tail estimate from Theorem~\ref{pt9} to the term on the right hand side of (\ref{pe23}). Thus for all $x\in K$ and for all $\epsilon\in(0,\epsilon_1]$ we have
\begin{eqnarray}
P^\epsilon(\mathcal{A}(\epsilon,\,x)^c)&\le&C_1\exp\left(\frac{-\delta^2}{\epsilon K_1}\right)\label{pe23a}\\
&\rightarrow&0\,\,\,\,\,\mbox{ as }\epsilon\rightarrow 0,\nonumber
\end{eqnarray}
where the numbers $C_1$ and $K_1$ in Peszat's exponential estimate (\ref{pe23a}) are positive real constants that do not depend on $\epsilon$ or on $x$.

\noindent We also have from Theorem~\ref{pt8}:
\begin{eqnarray*}
P^\epsilon(\mathcal{D}(\epsilon)^c)&\rightarrow&0 \,\,\,\,\,\mbox{ as }\epsilon\rightarrow 0.\,\,\,\,\,\,\,\,\qed
\end{eqnarray*}

\section{The upper bound}\label{ps6}

\noindent In this section we assume that (A3) holds and we prove Theorem~\ref{pt3} using the following proposition.

\begin{pro}\label{pp10}
Let $K\subset H$ be compact. Given $a>0$ and $\delta>0$ and $\phi\in L^2([0,1];U)$ there exists $\epsilon_0>0$ and $b>0$ such that for all $\epsilon\in(0,\epsilon_0]$ and for all $x\in K$ we have
\[
P\left\{\sup_{t\in[0,1]}|X^\epsilon_x(t)-z^\phi_x(t)|\geq\delta,\,\,\sup_{t\in[0,1]}\left|\epsilon^\frac{1}{2}W(t)-J\int_0^t\phi(s)\,ds\right|_{U_1}\le b\right\}\le\exp\left(-\frac{a}{\epsilon}\right).
\]
\end{pro}
The virtue of this proposition is that given positive $\delta$ the exponential bound on the right hand side has $a$, which we can choose to be as large as we please, in the numerator; the cost is the restriction on $\epsilon^\frac{1}{2}W$, but we have the large deviation principle of $\{\mathcal{L}(\epsilon^\frac{1}{2}W):\,\epsilon\in(0,1]\}$ to describe how these distributions behave. There is some work involved in arriving at the proof of Proposition~\ref{pp10} and this is left till the end. We only remark that we need several lemmas which use assumptions (A1), (A2) and (A3).

\medskip

\noindent \textbf{Proof of Theorem~\ref{pt3} assuming (A3).} This proof is almost identical to that of Peszat's Theorem 1.3 in~\cite{Pe}; we give all the details because of its importance.

\noindent Let $K$ be a compact subset of $H$. Fix $r>0$ and $\delta>0$ and $\gamma>0$. Let $a$ be a positive real number, to be specified later. By Proposition~\ref{pp10}, for each $\phi\in L^2([0,1];U)$ there exists $b_\phi>0$ and $\epsilon_\phi>0$ such that for all $\epsilon\in(0,\epsilon_\phi]$ and for all $x\in K$ we have
\begin{equation}\label{pe24}
P\left\{\sup_{t\in[0,1]}|X^\epsilon_x(t)-z^\phi_x(t)|\geq\delta,\,\sup_{t\in[0,1]}\left|\epsilon^\frac{1}{2}W(t)-J\int_0^t\phi(s)\,ds\right|_{U_1}\le b_\phi\right\}\le\exp\left(-\frac{a}{\epsilon}\right).
\end{equation}
Recall from equation (\ref{pe19}) that $\mathcal{I}_W$ is the rate function of the large deviation principle satisfied by $\{\mathcal{L}(\epsilon^\frac{1}{2}W):\,\epsilon\in(0,1]\}$. We have
\begin{eqnarray*}
\{\mathcal{I}_W\le r\}&=&\left\{u\in C([0,1];U_1):\,u(t)=J\int_0^t\psi(s)\,ds\,\,\,\mbox{ }\forall t\in[0,1],\right.\\
&&\,\,\,\,\,\mbox{\hspace{28mm}}\left.\mbox{ where }\psi\in L^2([0,1];U)\mbox{ and }\int_0^1|\psi(s)|_U^2\,ds\le 2r\right\}\\
&\subset&\bigcup_{\begin{array}{c}\scriptstyle{\psi\in L^2([0,1];U):}\\\scriptstyle{\int_0^1|\psi(s)|_U^2\,ds\le 2r}\end{array}}\left\{v\in C([0,1];U_1):\,\sup_{t\in[0,1]}\left|v(t)-J\int_0^t\psi(s)\,ds\right|_{U_1}<b_\psi\right\}.
\end{eqnarray*}
Since $\{\mathcal{I}_W\le r\}$ is a compact subset of $C([0,1];U_1)$, there exists a natural number $l$ and $\phi_1,\ldots,\phi_l\in L^2([0,1];U)$ such that $\int_0^1|\phi_j(s)|_U^2\,ds\le 2r$ for each $j\in\{1,\ldots,l\}$ and
\begin{equation}\label{pe25}
\{\mathcal{I}_W\le r\}\subset\bigcup_{j=1}^l\left\{v\in C([0,1];U_1): \sup_{t\in[0,1]}\left|v(t)-J\int_0^t\phi_j(s)\,ds\right|_{U_1}<b_{\phi_j}\right\}=:\mathcal{C}.
\end{equation}
For each $x\in H$ we may appeal to the definition of $\mathcal{C}$ in (\ref{pe25}) and write
\begin{eqnarray}
&&\!\!\!\!\!\!\!\!P\{X^\epsilon_x\notin B_{C([0,1];H)}(\{\mathcal{I}_x\le r\},\,\delta)\}\nonumber\\
&\le& P\{X^\epsilon_x\notin B_{C([0,1];H)}(\{\mathcal{I}_x\le r\},\,\delta),\,\epsilon^\frac{1}{2}W\in \mathcal{C}\}+P\{\epsilon^\frac{1}{2}W\notin\mathcal{C}\}\nonumber\\
&\le&\sum_{j=1}^lP\left\{X^\epsilon_x\notin B_{C([0,1];H)}(\{\mathcal{I}_x\le r\},\,\delta),\,\sup_{t\in[0,1]}\left|\epsilon^\frac{1}{2}W(t)-J\int_0^t\phi_j(s)\,ds\right|_{U_1}<b_{\phi_j}\right\}\nonumber\\
&&\mbox{\hspace{2mm}}{}+P\{\epsilon^\frac{1}{2}W\notin\mathcal{C}\}.\label{pe26}
\end{eqnarray}
Set $\epsilon_1:=\min\{\epsilon_{\phi_1},\ldots,\epsilon_{\phi_l}\}$. For each $j\in\{1,\ldots,l\}$ we have from inequality (\ref{pe24}) that for all $\epsilon\in(0,\epsilon_1]$ and for all $x\in K$
\begin{eqnarray}
&&\!\!\!\!\!\!\!\!\!\!\!\!P\left\{X^\epsilon_x\notin B_{C([0,1];H)}(\{\mathcal{I}_x\le r\},\,\delta),\,\sup_{t\in[0,1]}\left|\epsilon^\frac{1}{2}W(t)-J\int_0^t\phi_j(s)\,ds\right|_{U_1}< b_{\phi_j}\right\}\nonumber\\
&\le&\exp\left(-\frac{a}{\epsilon}\right).\label{pe27}
\end{eqnarray}
Since the open set $\mathcal{C}$ contains $\{\mathcal{I}_W\le r\}$, by the upper bound of the large deviation principle of the family $\{\mathcal{L}(\epsilon^\frac{1}{2}W):\,\epsilon\in(0,1]\}$ there exists $\epsilon_2>0$ such that for all $\epsilon\in(0,\epsilon_2]$
\begin{equation}\label{pe28}
P\{\epsilon^\frac{1}{2}W\notin\mathcal{C}\}\le\exp\left(\frac{-r+\frac{\gamma}{2}}{\epsilon}\right).
\end{equation}
Set $\epsilon_3:=\epsilon_1\wedge\epsilon_2$. Returning to inequality (\ref{pe26}), we have for all $x\in K$ and for all $\epsilon\in(0,\epsilon_3]$
\begin{eqnarray*}
P\{X^\epsilon_x\notin B_{C([0,1];H)}(\{\mathcal{I}_x\le r\},\,\delta)\}&\le&l\exp\left(-\frac{a}{\epsilon}\right)+\exp\left(\frac{-r+\frac{\gamma}{2}}{\epsilon}\right)\\
&\le&(l+1)\exp\left(\frac{-r+\frac{\gamma}{2}}{\epsilon}\right)
\end{eqnarray*}
when $a$ is taken as $r-\frac{\gamma}{2}$.

\noindent Finally set $\epsilon_0:=\epsilon_3\wedge\frac{\gamma}{2\ln(l+1)}$. $\qed$

\bigskip
\noindent Now we work towards proving Proposition~\ref{pp10}. In the following we make use of (A3):
\[
\Gamma:=\sup_{x\in H}\|G(x)\|_{L_2(U,H)}<\infty,
\]
as well as (A1) and (A2).

\noindent Fix $\phi\in L^2([0,1];U)$. For each $\epsilon\in(0,1]$ define $\tilde{F}_\epsilon:([0,1]\times H,\mathcal{B}_{[0,1]}\otimes\mathcal{B}_H)\rightarrow(H,\mathcal{B}_H)$ by
\begin{equation}\label{pe29}
\tilde{F}_\epsilon(s,x):=\epsilon F(\epsilon s,x)+G(x)\phi(s)\,\,\,\,\,\mbox{ }\forall s\in[0,1]\mbox{ and }\forall x\in H.
\end{equation}
It is not difficult to show that, for each $\epsilon\in(0,1]$, $\tilde{F}_\epsilon$ is measurable and
\begin{equation}\label{pe30}
|\tilde{F}_\epsilon(s,x)-\tilde{F}_\epsilon(s,y)|\le\nu_\epsilon(s)|x-y|\,\,\,\,\,\mbox{ }\forall x,y\in H\mbox{ and }\forall s\in[0,1]
\end{equation}
and
\begin{equation}\label{pe31}
|\tilde{F}_\epsilon(s,x)|\le\nu_\epsilon(s)(1+|x|)\,\,\,\,\,\mbox{ }\forall s\in[0,1]\mbox{ and }\forall x\in H,
\end{equation}
where $\nu_\epsilon(s):=\epsilon\nu(\epsilon s)+\Lambda|\phi(s)|_U$, $s\in[0,1]$, is a function in $L^2([0,1];\mathbb{R})$.

\smallskip

\noindent By Theorem~\ref{ptA1}, for each $\epsilon\in(0,1]$ and each $x\in H$ we may define $(Z^\epsilon_x(t))_{t\in[0,1]}$ as the continuous $(\mathcal{F}_t)$-predictable process such that
\begin{equation}\label{pe31a}
Z^\epsilon_x(t)=S(\epsilon t)x+\int_0^tS(\epsilon(t-s))\tilde{F}_\epsilon(s,Z^\epsilon_x(s))\,ds+\epsilon^\frac{1}{2}\int_0^tS(\epsilon(t-s))G(Z^\epsilon_x(s))\,dW(s)
\end{equation}
for all $t\in[0,1]$, $P$ a.e.. To prove Proposition~\ref{pp10} we will need some lemmas concerning the processes $(Z^\epsilon_x(t))_{t\in[0,1]}$. Lemmas~\ref{pl11} to~\ref{pl15} are ultimately used to prove Proposition~\ref{pp16}; Proposition~\ref{pp16}, in turn, is used to prove Proposition~\ref{pp10}. The one point in these lemmas where there is a difference from Peszat's methods worth remarking on is in Lemma~\ref{pl15}, where we use the equicontinuity assumption (A2) to control the size of one of the terms arising in the proof; Peszat's proof of the corresponding small noise asymptotics lemma~\cite[Lemma 6.4]{Pe} uses only continuity in the norm topology of $t\in(0,\infty)\mapsto S(t)$. For this reason, except for brief comments, we omit the proofs of Lemmas~\ref{pl11} to~\ref{pl14}.

\begin{lem}\label{pl11}
Given $a\in(0,\infty)$ and $R\in(0,\infty)$ there exists $D\in(0,\infty)$ such that for all $\epsilon\in(0,1]$ and for all $x\in B_H(0,R)$ we have
\[
P\{\sup_{t\in[0,1]}|Z^\epsilon_x(t)|\geq D\}\le \exp\left(-\frac{a}{\epsilon}\right).
\]
\end{lem}
To prove this lemma we firstly use Gronwall's Lemma; then, thanks to assumption (A3), we can use Theorem~\ref{pt9} to get an exponential bound.

\bigskip

\noindent We introduce some notation. Set
\[
t_{n,k}:=\frac{k}{2^n}\,\,\,\,\,\mbox{ for }n\in\mathbb{N}\mbox{ and }k=0,1,\ldots,2^n.
\]
We will be approximating $Z^\epsilon_x(t)$ using the values at the discrete times $t_{n,k}$.
\begin{lem}\label{pl12}
Given $a>0$ and $\delta>0$ there is a natural number $N$ such that for each $n\geq N$ there exists $\epsilon_n>0$ such that for all $\epsilon\in(0,\epsilon_n]$ and for all $x\in H$ we have
\[
P\left\{\sup_{k\in\{0,1,\ldots,2^n-1\}}\sup_{t\in[t_{n,k},t_{n,k+1}]}|\epsilon^\frac{1}{2}\int_{t_{n,k}}^tS(\epsilon(t-s))G(Z^\epsilon_x(s))\,dW(s)|\geq\delta\right\}\le\exp\left(-\frac{a}{\epsilon}\right).
\]
\end{lem}
The proof of this lemma uses Theorem~\ref{pt9} to get an exponential bound.

\bigskip

\noindent To simplify notation, for each natural number $n$ define the function
\[
\pi_n(t):=\left\{\begin{array}{cl}\frac{k}{2^n}&\mbox{ if }t\in\left(\frac{k}{2^n},\,\frac{k+1}{2^n}\right]\,\,,\,\,\,\,\,\mbox{ }k=0,1,\ldots,2^n-1\\
0&\mbox{ if }t=0.\end{array}\right.
\]

\noindent \begin{lem}\label{pl13}
Let $R>0$. Given $a>0$ and $\delta>0$ there is a natural number $n_0$ such that for each $n\geq n_0$ there exists $\epsilon_n>0$ such that for all $\epsilon\in(0,\epsilon_n]$ and for all $x\in B_H(0,R)$
\[
P\{\sup_{t\in[0,1]}|Z^\epsilon_x(t)-S(\epsilon(t-\pi_n(t)))Z^\epsilon_x(\pi_n(t))|\geq\delta\}\le\exp\left(-\frac{a}{\epsilon}\right).
\]
\end{lem}
This lemma is proved using Lemmas~\ref{pl11} and~\ref{pl12}.

\bigskip

\begin{lem}\label{pl14}
Let $R>0$. Given $a>0$ and $\delta>0$ there is a natural number $n_0$ such that for each $n\geq n_0$ there exists $\epsilon_n>0$ such that for all $\epsilon\in(0,\epsilon_n]$ and for all $x\in B_H(0,R)$ and for all $T\in[0,1]$ we have
\begin{eqnarray*}
&&\!\!\!\!\!\!\!P\left\{\epsilon^\frac{1}{2}\left|\int_0^TS(\epsilon(T-s))G(Z^\epsilon_x(s))\,dW(s)\right.\right.\mbox{\hspace{70mm}}\\
&&\mbox{\hspace{5mm}}{}\left.\left.-\int_0^TS(\epsilon(T-s))G(S(\epsilon(s-\pi_n(s)))Z^\epsilon_x(\pi_n(s)))\,dW(s)\right|\geq\delta\right\}\\
&\le&\exp\left(-\frac{a}{\epsilon}\right).
\end{eqnarray*}
\end{lem}
In the proof of Lemma~\ref{pl14} we use the right continuity of the filtration $(\mathcal{F}_t)$; it ensures that for given $x\in B_H(0,R)$, $n\in\mathbb{N}$ and $\epsilon\in(0,1]$ the random variable
\[
\tau_\rho(\omega):=\inf\{t\in[0,1]:\,|Z^\epsilon_x(t)(\omega)-S(\epsilon(t-\pi_n(t)))Z^\epsilon_x(\pi_n(t))(\omega)|\geq\rho\},\,\,\,\,\,\omega\in\Omega,
\]
where we set $\inf\emptyset=1$, is a $(\mathcal{F}_t)$-stopping time for any positive $\rho$. Introducing $\tau_\rho$ enables the use of Theorem~\ref{pt8} together with Lemma~\ref{pl13} in the proof.

\bigskip

\noindent\begin{lem}\label{pl15}
Given $a>0$ and $\delta>0$ and $0\le T_1<T_2\le 1$ and $R>0$ there exists $b>0$ and there exists $\epsilon_0\in(0,1]$ such that for each $x\in B_H(0,R)$ and for each $\epsilon\in(0,\epsilon_0]$
\begin{eqnarray*}
&&\!\!\!\!\!\!\!\!\!P\{\,|\epsilon^\frac{1}{2}\int_{T_1}^{T_2}S(\epsilon(T_2-s))G(S(\epsilon(s-T_1))Z^\epsilon_x(T_1))\,dW(s)|\geq\delta,\,\,\sup_{T_1\le t\le T_2}\epsilon^\frac{1}{2}|W(t)|_{U_1}\le b\,\}\\
&\le&\exp\left(\frac{-a}{\epsilon}\right).
\end{eqnarray*}
\end{lem}
\textbf{Proof.} Recall that $(g_k)$ is an orthonormal basis of $U$ and for each $n\in\mathbb{N}$ we define the projection in $U$:
\[
\Pi_n(u)=\sum_{k=1}^n \langle u,g_k\rangle_U g_k\,\,\,\,\,\mbox{ }\forall u\in U.
\]
In the course of this proof we choose numbers $D\in(0,\infty)$, $n\in\mathbb{N}$, $T_1<\tilde{T}_1<\tilde{T}_2<T_2$ and a partition $\tilde{T}_1=t_0<t_1<\cdots<t_l=\tilde{T}_2$ as well as $b\in(0,\infty)$ in order to control the size of the five terms on the right hand side of the inequality
\begin{eqnarray}
&&\!\!\!\!\!\!\!\!\!P\{\,\epsilon^\frac{1}{2}|\int_{T_1}^{T_2}S(\epsilon(T_2-s))G(S(\epsilon(s-T_1))Z^\epsilon_x(T_1))\,dW(s)|\geq\delta,\,\,\sup_{t\in[T_1,T_2]}\epsilon^\frac{1}{2}|W(t)|_{U_1}\le b\}\mbox{\hspace{1cm}}\nonumber\\
&\le&P\{|Z^\epsilon_x(T_1)|\geq D\}\nonumber\\
&&{}+P\{\,\epsilon^\frac{1}{2}|\int_{T_1}^{T_2}S(\epsilon(T_2-s))G(S(\epsilon(s-T_1))Z^\epsilon_x(T_1))(I_U-\Pi_n)\,dW(s)\,|\geq\frac{\delta}{4},\,|Z^\epsilon_x(T_1)|<D\}\nonumber\\
&&{}+P\left\{\,\epsilon^{\frac{1}{2}}\left|\int_{T_1}^{\tilde{T}_1}S(\epsilon(T_2-s))G(S(\epsilon(s-T_1))Z^\epsilon_x(T_1))\Pi_n\,dW(s)+{}\right.\right.\nonumber\\
&&\mbox{\hspace{20mm}}\left.\left.\int_{\tilde{T}_2}^{T_2}S(\epsilon(T_2-s))G(S(\epsilon(s-T_1))Z^\epsilon_x(T_1))\Pi_n\,dW(s)\right|\geq\frac{\delta}{4}\right\}\nonumber\\
&&{}+P\left\{\,\epsilon^\frac{1}{2}\left|\int_{\tilde{T}_1}^{\tilde{T}_2}S(\epsilon(T_2-s))G(S(\epsilon(s-T_1))Z^\epsilon_x(T_1))\Pi_n\,dW(s)\right.\right.\nonumber\\
&&\mbox{\hspace{20mm}}\left.{}-\int_{\tilde{T}_1}^{\tilde{T}_2}\sum_{j=0}^{l-1}1_{(t_j,t_{j+1}]}(s)S(\epsilon(T_2-t_j))G(S(\epsilon(t_j-T_1))Z^\epsilon_x(T_1))\Pi_n\,dW(s)\right|\geq\frac{\delta}{4},\nonumber\\
&&\mbox{\hspace{117mm}}\left.\begin{array}{c}\mbox{ }\\\mbox{ }\end{array}|Z^\epsilon_x(T_1)|<D\,\right\}\nonumber\\
&&{}+P\left\{\,\epsilon^\frac{1}{2}|\int_{\tilde{T}_1}^{\tilde{T}_2}\sum_{j=0}^{l-1}1_{(t_j,t_{j+1}]}(s)S(\epsilon(T_2-t_j))G(S(\epsilon(t_j-T_1))Z^\epsilon_x(T_1))\Pi_n\,dW(s)\,|\geq\frac{\delta}{4},\right.\nonumber\\
&&\mbox{\hspace{105mm}}\left.\sup_{t\in[T_1,T_2]}\epsilon^\frac{1}{2}|W(t)|_{U_1}\le b\,\right\}\nonumber\\
&=&\mbox{term 1}+\mbox{term 2}+\mbox{term 3}+\mbox{term 4}+\mbox{term 5}.\label{pe40a}
\end{eqnarray}
Let $\tilde{a}>a$.
\medskip

\noindent By Lemma~\ref{pl11} we can take $D\in(0,\infty)$ such that for all $x\in B_H(0,R)$ and for all $\epsilon\in(0,1]$
\begin{equation}\label{pe40b}
\mbox{term 1}:=P\{|Z^\epsilon_x(T_1)|\geq D\}\le\exp\left(-\frac{\tilde{a}}{\epsilon}\right).
\end{equation}

\medskip

\noindent Let $x\in B_H(0,R)$ and $\epsilon\in(0,1]$. Define the $(\mathcal{F}_t)$-stopping time
\begin{equation}\label{pe41}
\tau_{x,\epsilon}(\omega):=\left\{\begin{array}{ll}T_1&\mbox{ if }|Z^\epsilon_x(T_1)(\omega)|\geq D\\
1&\mbox{ otherwise.}\end{array}\right.
\end{equation}
We have
\begin{eqnarray}
&&\!\!\!\!\!\!\!\!\!\mbox{term 2\hspace{90mm}}\nonumber\\
&:=&\!P\{\,|\epsilon^\frac{1}{2}\int_{T_1}^{T_2}S(\epsilon(T_2-s))G(S(\epsilon(s-T_1))Z^\epsilon_x(T_1))(I_U-\Pi_n)\,dW(s)|\geq\frac{\delta}{4},\,|Z^\epsilon_x(T_1)|<D\,\}\nonumber\\
&\le&\!P\left\{\,\sup_{t\in[0,1]}|\epsilon^\frac{1}{2}\int_0^{t\wedge\tau_{x,\epsilon}}1_{(T_1,T_2]}(s)S(\epsilon(T_2-s))G(S(\epsilon(s-T_1))Z^\epsilon_x(T_1))(I_U-\Pi_n)\,dW(s)|\geq\frac{\delta}{4},\right.\nonumber\\
&&\mbox{\hspace{0.5cm}}\left.|Z^\epsilon_x(T_1)|<D\begin{array}{l}\mbox{ }\\\mbox{ }\end{array}\right\}.\,\,\,\,\,\,\,\,\,\,\,\,\label{pe42}
\end{eqnarray}
Since
\begin{eqnarray*}
&&\!\!\!\!\!\!\!\!\!\int_0^1\epsilon 1_{[0,\tau_{x,\epsilon}]}(s)1_{(T_1,T_2]}(s)\|S(\epsilon(T_2-s))G(S(\epsilon(s-T_1))Z^\epsilon_x(T_1))(I_U-\Pi_n)\|^2_{L_2(U,H)}\,ds\mbox{\hspace{2cm}}\\
&\le&\epsilon M^2\,\sup_{h\in B_H(0,MD)}\|G(h)(I_U-\Pi_n)\|^2_{L_2(U,H)}\,\,\,\,\,\mbox{ }P\mbox{ a.e.,}
\end{eqnarray*}
Theorem~\ref{pt8} yields an estimate of the term on the right hand side of inequality~(\ref{pe42}):
\begin{eqnarray*}
\mbox{term 2}
&\le&3\exp\left(-\frac{\delta^2}{64\epsilon M^2\,\sup_{h\in B_H(0,DM)}\|G(h)(I_U-\Pi_n)\|^2_{L_2(U,H)}}\right).
\end{eqnarray*}
By assumption (A1) we can now choose $n\in\mathbb{N}$ such that
\begin{equation}\label{pe40c}
\mbox{term 2}\le\exp\left(-\frac{\tilde{a}}{\epsilon}\right)\,\,\,\,\,\mbox{ }\forall x\in B_H(0,R)\mbox{ and }\forall\epsilon\in(0,1].
\end{equation}

\medskip

\noindent We choose $\tilde{T}_1$ and $\tilde{T}_2$ such that $T_1<\tilde{T}_1<\tilde{T}_2<T_2$.
Then again by Theorem~\ref{pt8} we have
\begin{eqnarray}
&&\!\!\!\!\!\!\!\!\!\mbox{term 3}\mbox{\hspace{9cm}}\nonumber\\
&:=&P\left\{\epsilon^\frac{1}{2}\left|\int_{T_1}^{\tilde{T}_1}S(\epsilon(T_2-s))G(S(\epsilon(s-T_1))Z^\epsilon_x(T_1))\Pi_n\,dW(s)\mbox{\hspace{20mm}}\right.\right.\nonumber\\
&&\mbox{\hspace{15mm}}\left.\left.{}+\int_{\tilde{T}_2}^{T_2}S(\epsilon(T_2-s))G(S(\epsilon(s-T_1))Z^\epsilon_x(T_1))\Pi_n\,dW(s)\right|\geq\frac{\delta}{4}\right\}\nonumber\\
&\le&P\{\sup_{t\in[0,1]}|\int_0^t1_{(T_1,\tilde{T}_1]\cup(\tilde{T}_2,T_2]}(s)\epsilon^\frac{1}{2}S(\epsilon(T_2-s))G(S(\epsilon(s-T_1))Z^\epsilon_x(T_1))\Pi_n\,dW(s)|\geq\frac{\delta}{4}\}\nonumber\\
&\le&3\exp\left(-\frac{\delta^2}{64\epsilon M^2\Gamma^2(\tilde{T}_1-T_1+T_2-\tilde{T}_2)}\right)\,,\,\mbox{ and for }\tilde{T}_1-T_1+T_2-\tilde{T}_2 \mbox{ small,}\nonumber\\
&\le&\exp\left(-\frac{\tilde{a}}{\epsilon}\right)\,\,\,\,\,\mbox{ }\forall x\in B_H(0,R)\mbox{ and }\forall\epsilon\in(0,1].\label{pe40d}
\end{eqnarray}

\medskip

\noindent Let $\mathcal{T}:=\{\tilde{T}_1=t_0<t_1<\cdots<t_l=\tilde{T}_2\}$ be a partition of $[\tilde{T}_1,\tilde{T}_2]$ and set $\Delta_{\mathcal{T}}:=\max\{t_{j+1}-t_j\,:\,j=0,1,\ldots,l-1\}$.

\noindent For $x\in B_H(0,R)$ and $\epsilon\in(0,1]$ define the $(\mathcal{F}_t)$-stopping time $\tau_{x,\epsilon}$ as in equation (\ref{pe41}). We have
\begin{eqnarray}
&&\!\!\!\!\!\!\!\!\!\mbox{term 4}\mbox{\hspace{7cm}}\nonumber\\
&:=&\!P\left\{\epsilon^\frac{1}{2}\left|\int_{\tilde{T}_1}^{\tilde{T}_2}S(\epsilon(T_2-s))G(S(\epsilon(s-T_1))Z^\epsilon_x(T_1))\Pi_n\,dW(s)\right.\right.\nonumber\\
&&\mbox{\hspace{10mm}}\left.{}-\int_{\tilde{T}_1}^{\tilde{T}_2}\sum_{j=0}^{l-1}1_{(t_j,t_{j+1}]}(s)S(\epsilon(T_2-t_j))G(S(\epsilon(t_j-T_1))Z^\epsilon_x(T_1))\Pi_n\,dW(s)\right|\geq\frac{\delta}{4},\nonumber\\
&&\mbox{\hspace{108mm}}\left.\begin{array}{c}\mbox{ }\\\mbox{ }\end{array}|Z^\epsilon_x(T_1)|<D\right\}\nonumber\\
&\le&\!P\left\{\sup_{t\in[0,1]}\epsilon^\frac{1}{2}\left|\int_0^{t\wedge\tau_{x,\epsilon}}\sum_{j=0}^{l-1}1_{(t_j,t_{j+1}]}(s)[S(\epsilon(T_2-s))G(S(\epsilon(s-T_1))Z^\epsilon_x(T_1))\right.\right.\nonumber\\
&&\mbox{\hspace{42mm}}\left.\left.{}\begin{array}{l}\mbox{ }\\\mbox{ }\\\mbox{}\end{array}-S(\epsilon(T_2-t_j))G(S(\epsilon(t_j-T_1))Z^\epsilon_x(T_1))]\Pi_n\,dW(s)\right|\geq\frac{\delta}{4}\right\}.\,\,\,\,\,\,\,\,\,\,\,\label{pe43}
\end{eqnarray}
In order to apply Theorem~\ref{pt8} to the right hand side of (\ref{pe43}) we observe that
\begin{eqnarray*}
&&\!\!\!\!\!\!\!\!\!\int_0^1\sum_{j=0}^{l-1}1_{(t_j,t_{j+1}]}(s)1_{[0,\tau_{x,\epsilon}]}(s)\|[S(\epsilon(T_2-s))G(S(\epsilon(s-T_1))Z^\epsilon_x(T_1))\\
&&\mbox{\hspace{39mm}}{}-S(\epsilon(T_2-t_j))G(S(\epsilon(t_j-T_1))Z^\epsilon_x(T_1))]\Pi_n\|^2_{L_2(U,H)}\,ds\\
&\le&2(M^2 \Lambda^2 D^2+\Gamma^2)(\zeta(\Delta_\mathcal{T}))^2,
\end{eqnarray*}
where
\begin{eqnarray*}
\zeta(\Delta_\mathcal{T})&:=&\sup\left\{\|S(\eta r)-S(\eta s)\|_{L(H,H)}:r,s\in[(\tilde{T}_1-T_1)\wedge(T_2-\tilde{T}_2),1]\begin{array}{c}\mbox{ }\\\mbox{ }\end{array}\right.\mbox{\hspace{3cm}}\\
&&\left.\begin{array}{c}\mbox{ }\\\mbox{ }\end{array}\mbox{\hspace{43mm} and }|r-s|\le\Delta_{\mathcal{T}}\mbox{ and }\eta\in(0,1]\right\}.
\end{eqnarray*}
The number $\zeta(\Delta_\mathcal{T})$ goes to $0$ as $\Delta_\mathcal{T}\rightarrow 0$ since, by (A2), the family of functions
\[
\{t\in[(\tilde{T}_1-T_1)\wedge(T_2-\tilde{T}_2),1]\mapsto S(\eta t)\in L(H,H),\,\,\,\mbox{ }\eta\in(0,1]\}
\]
is uniformly equicontinuous in the norm topology.

\noindent We now choose partition $\mathcal{T}=\{\tilde{T}_1=t_0<t_1<\cdots<t_l=\tilde{T}_2\}$ such that $\Delta_\mathcal{T}$ satisfies
\[
\ln 3-\frac{\delta^2}{128(M^2\Lambda^2D^2+\Gamma^2)(\zeta(\Delta_\mathcal{T}))^2}\le\,\,\,-\tilde{a}.
\]
Then from inequality (\ref{pe43}) and Theorem~\ref{pt8} we have
\begin{equation}\label{pe40e}
\mbox{term 4}\le\exp\left(-\frac{\tilde{a}}{\epsilon}\right)\,\,\,\,\,\mbox{ }\forall x\in B_H(0,R)\mbox{ and }\forall\epsilon\in(0,1].
\end{equation}

\medskip

\noindent Finally we consider term 5. Recall that by definition of the inner product $\langle\cdot,\cdot\rangle_{U_1}$ in $U_1$, the bounded linear operator from $U_1$ into $U$
\[
\Pi^1_nu:=\sum_{k=1}^n\langle u,\lambda_k^{-2}Jg_k\rangle_{U_1}\,g_k\,\,,\,\,\,\mbox{ }u\in U_1,
\]
satisfies $\Pi_n^1\,Ju=\Pi_nu\,\,\,\mbox{ }\forall u\in U$. We will use the result
\begin{equation}\label{pe44a}
\int_0^11_{(c,d]}(s)\Phi\circ J\,dW(s)=\Phi(W(d)-W(c))\,\,\,\,\,\mbox{ }P\mbox{ a.e.}
\end{equation}
when $0\le c<d\le 1$ and $\Phi:(\Omega,\mathcal{F}_c)\rightarrow (L_2(U_1,H),\mathcal{B}_{L_2(U_1,H)})$ is $\mathcal{F}_c$ measurable and $E[\|\Phi\|^2_{L_2(U_1,H)}]<\infty$.

\noindent We have for each $\epsilon\in(0,1]$ and $x\in H$
\begin{eqnarray}
&&\!\!\!\!\!\!\!\epsilon^\frac{1}{2}|\int_{\tilde{T}_1}^{\tilde{T}_2}\sum_{j=0}^{l-1}1_{(t_j,t_{j+1}]}(s)S(\epsilon(T_2-t_j))G(S(\epsilon(t_j-T_1))Z^\epsilon_x(T_1))\Pi_n\,dW(s)|\mbox{\hspace{20mm}}\nonumber\\
&=&\epsilon^\frac{1}{2}|\sum_{j=0}^{l-1}S(\epsilon(T_2-t_j))G(S(\epsilon(t_j-T_1))Z^\epsilon_x(T_1))\Pi_n^1(W(t_{j+1})-W(t_j))|\,\,\,\mbox{ }P\mbox{ a.e.}\,\,\,\,\,\label{pe44b}\\
&\le&2lM\Gamma\|\Pi_n^1\|_{L(U_1,U)}\epsilon^\frac{1}{2}\sup_{t\in[T_1,T_2]}|W(t)|_{U_1}.\label{pe45}
\end{eqnarray}
Equality (\ref{pe44b}) follows from equality (\ref{pe44a}). We choose $0<b<\frac{\delta}{8lM\Gamma\|\Pi^1_n\|_{L(U_1,U)}}$, then for each $\epsilon\in(0,1]$ and each $x\in B_H(0,R)$
\begin{eqnarray}
\mbox{term 5}
&=&P\left\{\epsilon^\frac{1}{2}|\int_{\tilde{T}_1}^{\tilde{T}_2}\sum_{j=0}^{l-1}1_{(t_j,t_{j+1}]}(s)S(\epsilon(T_2-t_j))G(S(\epsilon(t_j-T_1))Z^\epsilon_x(T_1))\Pi_n\,dW(s)|\geq\frac{\delta}{4},\right.\nonumber\\
&&\mbox{\hspace{85mm}}\left.\sup_{t\in[T_1,T_2]}\epsilon^\frac{1}{2}|W(t)|_{U_1}\le b\right\}=0,\nonumber\\
&&\mbox{\hspace{105mm}}\label{pe40f}
\end{eqnarray}
by inequality (\ref{pe45}).

\medskip

\noindent With $b$ chosen as in the last paragraph, we combine inequalities (\ref{pe40a}), (\ref{pe40b}), (\ref{pe40c}), (\ref{pe40d}), (\ref{pe40e}) and (\ref{pe40f}) to obtain the result of the lemma. $\qed$

\bigskip

\noindent\begin{pro}\label{pp16}
Let $R\in(0,\infty)$. Given $a>0$ and $\delta>0$ there exist $b>0$ and $\epsilon_0>0$ such that for all $x\in B_H(0,R)$ and for all $\epsilon\in(0,\epsilon_0]$ we have
\[
P\left\{\sup_{t\in[0,1]}\left|\epsilon^\frac{1}{2}\int_0^tS(\epsilon(t-s))G(Z^\epsilon_x(s))\,dW(s)\right|\geq\delta,\,\,\,\mbox{ }\sup_{t\in[0,1]}|\epsilon^\frac{1}{2}W(t)|_{U_1}\le b\right\}\le\exp\left(-\frac{a}{\epsilon}\right).
\]
\end{pro}
\textbf{Proof.} Let $\tilde{a}>a$. Let $n$ be a natural number. For each $k\in\{0,1,\ldots,2^n-1\}$ and $t\in[t_{n,k},t_{n,k+1}]$ we have
\begin{eqnarray}
\left|\epsilon^\frac{1}{2}\int_0^tS(\epsilon(t-s))G(Z^\epsilon_x(s))\,dW(s)\right|&\le&\left|\epsilon^\frac{1}{2}S(\epsilon(t-t_{n,k}))\int_0^{t_{n,k}}S(\epsilon(t_{n,k}-s))G(Z^\epsilon_x(s))\,dW(s)\right|\nonumber\\
&&\mbox{\hspace{5mm}}{}+\left|\epsilon^\frac{1}{2}\int_{t_{n,k}}^tS(\epsilon(t-s))G(Z^\epsilon_x(s))\,dW(s)\right|\,\,\,\mbox{ }P\mbox{ a.e..}\mbox{\hspace{10mm}}\label{pe46}
\end{eqnarray}
It follows that
\begin{eqnarray*}
&&\!\!\!\!\!\!\!\!\!\sup_{t\in[0,1]}\left|\epsilon^\frac{1}{2}\int_0^tS(\epsilon(t-s))G(Z^\epsilon_x(s))\,dW(s)\right|\mbox{\hspace{60mm}}\\
&\le&M\epsilon^\frac{1}{2}\,\sup_{0<k\le 2^n-1}\left|\int_0^{t_{n,k}}S(\epsilon(t_{n,k}-s))G(Z^\epsilon_x(s))\,dW(s)\right.\\
&&\left.\mbox{\hspace{25mm}}{}-\int_0^{t_{n,k}}S(\epsilon(t_{n,k}-s))G(S(\epsilon(s-\pi_n(s)))Z^\epsilon_x(\pi_n(s)))\,dW(s)\right|\\
&&\mbox{\hspace{2mm}}{}+M\epsilon^\frac{1}{2}\,\sup_{0<k\le 2^n-1}\left|\int_0^{t_{n,k}}S(\epsilon(t_{n,k}-s))G(S(\epsilon(s-\pi_n(s)))Z^\epsilon_x(\pi_n(s)))\,dW(s)\right|\\
&&\mbox{\hspace{2mm}}{}+\epsilon^\frac{1}{2}\,\sup_{0\le k\le 2^n-1}\,\sup_{t\in[t_{n,k},t_{n,k+1}]}\left|\int_{t_{n,k}}^tS(\epsilon(t-s))G(Z^\epsilon_x(s))\,dW(s)\right|\,\,\,\mbox{ }P\mbox{ a.e..}
\end{eqnarray*}
By Lemma~\ref{pl12} and Lemma~\ref{pl14} respectively, there exist a natural number $n_0$ and a positive number $\epsilon_0$ such that
\begin{enumerate}
\item for all $x\in H$ and for all $\epsilon\in(0,\epsilon_0]$ we have
\[
P\left\{\sup_{0\le k\le 2^{n_0}-1}\,\sup_{t\in[t_{n_0,k},t_{n_0,k+1}]}\left|\epsilon^\frac{1}{2}\int_{t_{n_0,k}}^tS(\epsilon(t-s))G(Z^\epsilon_x(s))\,dW(s)\right|\geq\frac{\delta}{3}\right\}\le\exp\left(-\frac{\tilde{a}}{\epsilon}\right)
\]
\item and for each $k\in\{1,\ldots,2^{n_0}-1\}$ and for all $x\in B_H(0,R)$ and for all $\epsilon\in(0,\epsilon_0]$ we have
\begin{eqnarray*}
&&\!\!\!\!\!P\left\{\epsilon^\frac{1}{2}\left|\int_0^{t_{n_0,k}}S(\epsilon(t_{n_0,k}-s))G(Z^\epsilon_x(s))\,dW(s)\right.\right.\mbox{\hspace{60mm}}\\
&&\mbox{\hspace{2mm}}{}\left.\left.-\int_0^{t_{n_0,k}}S(\epsilon(t_{n_0,k}-s))G(S(\epsilon(s-\pi_{n_0}(s)))Z^\epsilon_x(\pi_{n_0}(s)))\,dW(s)\right|\geq\frac{\delta}{3M}\right\}\le\exp\left(-\frac{\tilde{a}}{\epsilon}\right).
\end{eqnarray*}
\end{enumerate}
Hence for arbitrary $b>0$ and for all $x\in B_H(0,R)$ and for all $\epsilon\in(0,\epsilon_0]$ we have
\begin{eqnarray}
&&\!\!\!\!\!\!\!\!\!P\left\{\sup_{t\in[0,1]}\left|\epsilon^\frac{1}{2}\int_0^tS(\epsilon(t-s))G(Z^\epsilon_x(s))\,dW(s)\right|\geq\delta,\,\,\,\sup_{t\in[0,1]}|\epsilon^\frac{1}{2}W(t)|_{U_1}\le b\right\}\mbox{\hspace{20mm}}\nonumber\\
&\le&P\left\{\epsilon^\frac{1}{2}\sup_{0<k\le 2^{n_0}-1}\left|\int_0^{t_{n_0,k}}S(\epsilon(t_{n_0,k}-s))G(Z^\epsilon_x(s))\,dW(s)\right.\right.\nonumber\\
&&\mbox{\hspace{27mm}}\left.\left.{}-\int_0^{t_{n_0,k}}S(\epsilon(t_{n_0,k}-s))G(S(\epsilon(s-\pi_{n_0}(s)))Z^\epsilon_x(\pi_{n_0}(s)))\,dW(s)\right|\geq\frac{\delta}{3M}\right\}\nonumber\\
&&\mbox{\hspace{0mm}}{}+P\left\{\epsilon^\frac{1}{2}\sup_{0\le k\le 2^{n_0}-1}\,\sup_{t\in[t_{n_0,k},t_{n_0,k+1}]}\left|\int_{t_{n_0,k}}^tS(\epsilon(t-s))G(Z^\epsilon_x(s))\,dW(s)\right|\geq\frac{\delta}{3}\right\}\nonumber\\
&&\mbox{\hspace{0mm}}{}+P\left\{\epsilon^\frac{1}{2}\sup_{0<k\le 2^{n_0}-1}\left|\int_0^{t_{n_0,k}}S(\epsilon(t_{n_0,k}-s))G(S(\epsilon(s-\pi_{n_0}(s)))Z^\epsilon_x(\pi_{n_0}(s)))\,dW(s)\right|\geq\frac{\delta}{3M},\right.\nonumber\\
&&\mbox{\hspace{105mm}}\left.\sup_{t\in[0,1]}|\epsilon^\frac{1}{2}W(t)|_{U_1}\le b\right\}\nonumber\\
&\le&2^{n_0}\exp\left(-\frac{\tilde{a}}{\epsilon}\right)\nonumber\\
&&\mbox{\hspace{0mm}}{}+\sum_{j=0}^{2^{n_0}-2}P\left\{\epsilon^\frac{1}{2}\left|\int_{t_{n_0,j}}^{t_{n_0,j+1}}S(\epsilon(t_{n_0,j+1}-s))G(S(\epsilon(s-t_{n_0,j}))Z^\epsilon_x(t_{n_0,j}))\,dW(s)\right|\geq\frac{\delta}{3M^2(2^{n_0}-1)},\right.\nonumber\\
&&\mbox{\hspace{105mm}}\left.\sup_{t\in[0,1]}|\epsilon^\frac{1}{2}W(t)|_{U_1}\le b\right\}.\,\,\,\,\,\,\,\,\,\,\,\,\label{pe47}
\end{eqnarray}
Now we can use Lemma~\ref{pl15} to choose a suitable $b$ to put in inequality~(\ref{pe47}), which completes the proof. $\qed$

\bigskip

\noindent We can now prove Proposition~\ref{pp10}.

\noindent\textbf{Proof of Proposition~\ref{pp10}.} Let $K\subset H$ be compact. Fix $a>0$ and $\delta>0$ and $\phi\in L^2([0,1];U)$. For $\epsilon\in(0,1]$ and $x\in K$ and $b$ a positive real number which will be specified later, we set
\[
\mathcal{D}(\epsilon,x,b):=\left\{\,\sup_{t\in[0,1]}|X^\epsilon_x(t)-z^\phi_x(t)|\geq\delta,\,\,\,\sup_{t\in[0,1]}\left|\epsilon^\frac{1}{2}W(t)-J\int_0^t\phi(s)\,ds\right|_{U_1}\le b\right\}.
\]
As in equation~(\ref{pe20}), define the process
\[
W^\epsilon(t):= W(t)-\epsilon^{-\frac{1}{2}}J\int_0^t\phi(s)\,ds\,\,\,\,\,\mbox{ }\forall t\in[0,1].
\]
By~\cite[Theorem 10.14]{DPZ}, $(W^\epsilon(t))_{t\in[0,1]}$ is a Wiener process with respect to filtration $(\mathcal{F}_t)$ on probability space $(\Omega,\mathcal{F},P^\epsilon)$ where
\[
dP^\epsilon(\omega)=\exp\left(\epsilon^{-\frac{1}{2}}\int_0^1\langle\phi(s),\cdot\rangle_U\,dW(s)-\frac{1}{2\epsilon}\int_0^1|\phi(s)|_U^2\,ds\right)\,dP(\omega)
\]
and $P^\epsilon(W^\epsilon(1))^{-1}=P(W(1))^{-1}$.

\noindent For $\lambda>0$ set
\[
\mathcal{M}(\epsilon,\lambda):=\left\{\omega\in\Omega:\int_0^1\langle\phi(s),\cdot\rangle_U\,dW(s)(\omega)\geq-\frac{\lambda}{\epsilon^\frac{1}{2}}\right\}.
\]
We have
\begin{equation}\label{pe48}
P(\mathcal{D}(\epsilon,x,b))\le P(\mathcal{D}(\epsilon,x,b)\cap\mathcal{M}(\epsilon,\lambda))+P(\mathcal{M}(\epsilon,\lambda)^c).
\end{equation}
By Theorem~\ref{pt8} we have immediately
\begin{eqnarray}
P(\mathcal{M}(\epsilon,\lambda)^c)&\le&P\left\{\sup_{t\in[0,1]}\left|\int_0^t\langle\phi(s),\cdot\rangle_U\,dW(s)\right|>\frac{\lambda}{\epsilon^\frac{1}{2}}\right\}\nonumber\\
&\le&3\exp\left(-\frac{\lambda^2}{4\epsilon\int_0^1|\phi(s)|^2_U\,ds}\right).\label{pe49}
\end{eqnarray}
The rest of the proof involves finding an exponential bound for the first term on the right hand side of~(\ref{pe48}). We have
\begin{eqnarray}
&&\!\!\!\!\!\!\!\!\!P(\mathcal{D}(\epsilon,x,b)\cap\mathcal{M}(\epsilon,\lambda))\mbox{\hspace{70mm}}\nonumber\\
&=&\int_\Omega 1_{\mathcal{D}(\epsilon,x,b)\cap\mathcal{M}(\epsilon,\lambda)}(\omega)\exp\left(-\epsilon^{-\frac{1}{2}}\int_0^1\langle\phi(s),\cdot\rangle_U\,dW(s)(\omega)+\frac{1}{2\epsilon}\int_0^1|\phi(s)|^2_U\,ds\right)\,dP^\epsilon(\omega)\nonumber\\
&\le&\exp\left(\frac{\lambda}{\epsilon}+\frac{1}{2\epsilon}\int_0^1|\phi(s)|_U^2\,ds\right)\,P^\epsilon\left\{\sup_{t\in[0,1]}|X^\epsilon_x(t)-z^\phi_x(t)|\geq\delta,\,\,\,\sup_{t\in[0,1]}\epsilon^\frac{1}{2}|W^\epsilon(t)|_{U_1}\le b\right\}.\mbox{\hspace{11mm}}\label{pe50}
\end{eqnarray}
By Lemma~\ref{pAl1} we have
\begin{eqnarray*}
X^\epsilon_x(t)&=&S(\epsilon t)x+\epsilon\int_0^tS(\epsilon(t-s))F(\epsilon s, X^\epsilon_x(s))\,ds+\epsilon^\frac{1}{2}\int_0^tS(\epsilon(t-s))G(X^\epsilon_x(s))\,dW^\epsilon(s)\\
&&\mbox{\hspace{14mm}}{}+\int_0^tS(\epsilon(t-s))G(X^\epsilon_x(s))\phi(s)\,ds\,\,\,\,\,\mbox{\hspace{25mm}}\forall t\in[0,1]\,\,\,\mbox{ }P\mbox{ a.e..}
\end{eqnarray*}
Recall that in equation (\ref{pe29}) we defined
\[
\tilde{F}_\epsilon(s,y):=\epsilon F(\epsilon s,y)+G(y)\phi(s)\,\,\,\,\,\mbox{ }\forall(s,y)\in [0,1]\times H.
\]
Thus $(X^\epsilon_x(t):(\Omega,\mathcal{F}_t,P^\epsilon)\rightarrow (H,\mathcal{B}_H))_{t\in[0,1]}$ is the solution of the equation
\begin{equation}\label{pe51a}
\left\{\begin{array}{rcl}dX^\epsilon(t)&=&(\epsilon AX^\epsilon(t)+\tilde{F}_\epsilon(t,X^\epsilon(t)))\,dt+\epsilon^\frac{1}{2}G(X^\epsilon(t))\,dW^\epsilon(t)\,\,\,\,\,\mbox{ }t\in(0,1]\\
X^\epsilon(0)&=&x.\end{array}\right.
\end{equation}
Let $(Z^\epsilon_x(t):(\Omega,\mathcal{F}_t,P)\rightarrow (H,\mathcal{B}_H))_{t\in[0,1]}$ be the solution of the equation
\begin{equation}\label{pe51b}
\left\{\begin{array}{rcl}dZ^\epsilon(t)&=&(\epsilon AZ^\epsilon(t)+\tilde{F}_\epsilon(t,Z^\epsilon(t)))\,dt+\epsilon^\frac{1}{2}G(Z^\epsilon(t))\,dW(t)\,\,\,\,\,\mbox{ }t\in(0,1]\\
Z^\epsilon(0)&=&x.\end{array}\right.
\end{equation}
By Proposition~\ref{ppA1} we have the equality of the distributions on $(C([0,1];H\oplus U_1),\mathcal{B}_{C([0,1];H\oplus U_1)})$:
\[
P^\epsilon(X^\epsilon_x,W^\epsilon)^{-1}=P(Z^\epsilon_x,W)^{-1};
\]
here trajectory-valued random variables are defined as in equation (\ref{pe55}). Thus we have
\begin{eqnarray}
&&\!\!\!\!\!\!\!\!\!P^\epsilon\left\{\sup_{t\in[0,1]}|X^\epsilon_x(t)-z^\phi_x(t)|\geq\delta ,\,\,\,\sup_{t\in[0,1]}\epsilon^\frac{1}{2}|W^\epsilon(t)|_{U_1}\le b\right\}\mbox{\hspace{20mm}}\nonumber\\
&=&P\left\{\sup_{t\in[0,1]}|Z^\epsilon_x(t)-z^\phi_x(t)|\geq\delta ,\,\,\,\sup_{t\in[0,1]}\epsilon^\frac{1}{2}|W(t)|_{U_1}\le b\right\}.\label{pe52}
\end{eqnarray}
We have the inequality
\begin{eqnarray*}
&&\!\!\!\!\!\!\!\!\!\sup_{t\in[0,1]}|Z^\epsilon_x(t)-z^\phi_x(t)|^2\mbox{\hspace{90mm}}\\
&\le&6\left[\sup_{r\in[0,\epsilon]}|S(r)x-x|^2+\epsilon M^2(\textstyle{\int_0^1\nu^2(s)\,ds})\int_0^1(1+|z^\phi_x(s)|)^2\,ds\right.\\
&&\mbox{\hspace{2mm}}{}+\sup\{\|(S(r)-I_H)G(z^\phi_x(s))\|_{L_2(U,H)}:\,s\in[0,1]\mbox{ and }r\in[0,\epsilon]\}^2\int_0^1|\phi(s)|^2_U\,ds\\
&&\mbox{\hspace{2mm}}\left.{}+\sup_{r\in[0,1]}\left|\epsilon^\frac{1}{2}\int_0^rS(\epsilon(r-s))G(Z^\epsilon_x(s))\,dW(s)\right|^2\right]\\
&&\mbox{\hspace{40mm}}{}\times\exp\left(6M^2\left(\textstyle{\int_0^1\nu^2(s)\,ds}+\Lambda^2\textstyle{\int_0^1|\phi(s)|_U^2\,ds}\right)\right)\,\,\,\,\,P\mbox{ a.e.,}
\end{eqnarray*}
which is obtained in the same way as inequality~(\ref{e24a}). Then, since $K$ is compact, there exists $\epsilon_1>0$ such that for all $x\in K$ and for all $\epsilon\in(0,\epsilon_1]$ we have
\begin{eqnarray*}
&&\!\!\!\!\!\!\!\!\!P\left\{\sup_{t\in[0,1]}|Z^\epsilon_x(t)-z^\phi_x(t)|\geq\delta,\,\,\sup_{t\in[0,1]}\epsilon^\frac{1}{2}|W(t)|_{U_1}\le b\right\}\mbox{\hspace{10mm}}\\
&\le&P\left\{\sup_{r\in[0,1]}\left|\epsilon^\frac{1}{2}\int_0^rS(\epsilon(r-s))G(Z^\epsilon_x(s))\,dW(s)\right|\geq\frac{\delta}{c}\,,\,\,\,\sup_{t\in[0,1]}\epsilon^\frac{1}{2}|W(t)|_{U_1}\le b\right\},
\end{eqnarray*}
where $c:=3\exp(3M^2(\textstyle{\int_0^1\nu^2(s)\,ds}+\Lambda^2\int_0^1|\phi(s)|_U^2\,ds))$.

\noindent Given $\tilde{a}>a$, by Proposition~\ref{pp16} there exist positive real numbers $b$ and $\epsilon_2$ such that for all $\epsilon\in(0,\epsilon_2]$ and for all $x\in K$ we have
\[
P\left\{\sup_{t\in[0,1]}\left|\epsilon^\frac{1}{2}\int_0^tS(\epsilon(t-s))G(Z^\epsilon_x(s))\,dW(s)\right|\geq\frac{\delta}{c}\,,\,\,\,\sup_{t\in[0,1]}|\epsilon^\frac{1}{2}W(t)|_{U_1}\le b\right\}\le\exp\left(-\frac{\tilde{a}}{\epsilon}\right).
\]
Thus for all $x\in K$ and for all $\epsilon\in(0,\epsilon_1\wedge\epsilon_2]$ we have
\begin{equation}\label{pe53}
P\left\{\sup_{t\in[0,1]}|Z^\epsilon_x(t)-z^\phi_x(t)|\geq\delta\,,\,\,\,\sup_{t\in[0,1]}|\epsilon^\frac{1}{2}W(t)|_{U_1}\le b\right\}\le\exp\left(-\frac{\tilde{a}}{\epsilon}\right).
\end{equation}
Now from inequalities~(\ref{pe50}) and~(\ref{pe52}) and~(\ref{pe53}) we have for all $x\in K$ and for all $\epsilon\in(0,\epsilon_1\wedge\epsilon_2]$:
\begin{equation}\label{pe54}
P(\mathcal{D}(\epsilon,x,b)\cap\mathcal{M}(\epsilon,\lambda))\le\exp\left(\frac{\lambda+\frac{1}{2}\int_0^1|\phi(s)|_U^2\,ds-\tilde{a}}{\epsilon}\right).
\end{equation}
By firstly choosing $\lambda$ such that $\frac{-\lambda^2}{4\int_0^1|\phi(s)|_U^2\,ds}+\ln 3<-a$ and then choosing $\tilde{a}$ such that $\lambda+\frac{1}{2}\int_0^1|\phi(s)|^2_U\,ds-\tilde{a}<-a$, we see, on combining inequalities~(\ref{pe48}), (\ref{pe49}) and (\ref{pe54}), that the proof is complete. $\qed$

\section{Concluding remarks}\label{ps7}
By modifying Peszat's methods for the small noise asymptotics problem, we have obtained a large deviation principle describing the small time asymptotics of the continuous $H$-valued trajectories of the solution of equation~(\ref{pe4}). Peszat actually considered the situation where a Banach space $E$ is embedded in $H$ and $S(t)(H)\subset E$ for each $t>0$ and the restriction of $(S(t))_{t\geq 0}$ to $E$ gives a strongly continuous semigroup on $E$. When the Banach space $E$ is different from $H$, Peszat's methods for small noise asymptotics in the space of continuous $E$-valued functions, $C([0,1];E)$, cannot easily be adapted to give small time asymptotics results in $C([0,1];E)$. For example, when working in $C([0,1];E)$, $\kappa$ in Theorem~\ref{pt9} is defined with the operator norm $\|\cdot\|_{L(H,E)}$ in place of $\|\cdot\|_{L(H,H)}$ and, since we consider stochastic convolutions of semigroups $(S(\epsilon t))_{t\geq 0}$ depending on $\epsilon$, $\kappa$ may depend on $\epsilon$ in a troublesome way.

\section{Appendix}\label{psa}

Some results which we used in the main text are placed here to avoid cluttering the path leading to the proofs of the theorems.

\medskip

\noindent Let Hilbert spaces $H$, $U$ and $U_1$ be as defined in Section~\ref{ps1}.

\medskip

\noindent Let $T\in(0,\infty)$. Let $(\Omega,\mathcal{F},P)$ be a probability space. Let $(\mathcal{F}_t)_{t\geq 0}$ be a filtration of sub $\sigma$-algebras of $\mathcal{F}$ such that all sets in $\mathcal{F}$ of $P$ measure zero are in $\mathcal{F}_0$ and let $(W(t):(\Omega,\mathcal{F}_t,P)\rightarrow (U_1,\mathcal{B}_{U_1}))_{t\geq 0}$ be a $U_1$-valued $(\mathcal{F}_t)$-Wiener process such that $\mathcal{L}(W(1))$ has reproducing kernel Hilbert space $U$. Let $(S(t))_{t\geq 0}$ be a strongly continuous semigroup of bounded linear operators on $H$. Let the function $F:([0,T]\times H,\mathcal{B}_{[0,T]}\otimes\mathcal{B}_{H})\rightarrow (H,\mathcal{B}_H)$ be measurable and suppose there is a function $\theta\in L^2([0,T];\mathbb{R})$ such that
\begin{eqnarray}
|F(t,x)-F(t,y)|&\le &\theta(t)|x-y|\,\,\,\mbox{ }\forall t\in[0,T]\mbox{ and }\forall x,y\in H\mbox{ and }\label{pe6a}\\
|F(t,x)|&\le&\theta(t)(1+|x|)\,\,\,\mbox{ }\forall t\in[0,T]\mbox{ and }\forall x\in H.\label{pe6b}
\end{eqnarray}
Let $G:H\rightarrow L_2(U,H)$ be Lipschitz continuous. Let $x\in H$.
\begin{thm}[Existence, uniqueness and continuity of solutions]\label{ptA1}
There exists a $(\mathcal{F}_t)$-predictable process $(X(t))_{t\in[0,T]}$, unique up to equivalence among processes satisfying
\[
P\{\int_0^T|X(t)|^2\,dt<\infty\}=1,
\]
such that
\begin{equation}\label{e3.73aa}
X(t)=S(t)x+\int_0^tS(t-s)F(s,X(s))\,ds+\int_0^tS(t-s)G(X(s))\,dW(s)\,\,\,\mbox{ }P\mbox{ a.e.}
\end{equation}
for each $t\in[0,T]$. Moreover it has a continuous version and $\sup_{t\in[0,T]}E[|X(t)|^p]$ is finite for each $p\in[2,\infty)$.
\end{thm}
\smallskip

\noindent The proof of this theorem is omitted as it is almost identical to the proof of~\cite[Theorem 7.4]{DPZ}.

\bigskip

\noindent Let $P^\prime$ be another probability measure on the $\sigma$-algebra $\mathcal{F}$ of subsets of $\Omega$. Let $(\mathcal{G}_t)_{t\geq 0}$ be a filtration of sub $\sigma$-algebras of $\mathcal{F}$ such that all sets in $\mathcal{F}$ of $P^\prime$ measure zero are in $\mathcal{G}_0$. Let $(V(t):(\Omega,\mathcal{G}_t,P^\prime)\rightarrow (U_1,\mathcal{B}_{U_1}))_{t\geq 0}$ be a $U_1$-valued $(\mathcal{G}_t)$-Wiener process such that $\mathcal{L}(V(1))=\mathcal{L}(W(1))$.

\smallskip

\noindent Suppose that $(X(t))_{t\in[0,T]}$ is the continuous $(\mathcal{F}_t)$-predictable process in Theorem~\ref{ptA1} and suppose that $(Y(t))_{t\in[0,T]}$ is the continuous $(\mathcal{G}_t)$-predictable process satisfying
\begin{equation}\label{e3.73bb}
Y(t)=S(t)x+\int_0^tS(t-s)F(s,Y(s))\,ds+\int_0^tS(t-s)G(Y(s))\,dV(s)\,\,\,\mbox{ }P^\prime\mbox{ a.e.}
\end{equation}
for each $t\in[0,T]$. Let $H\oplus U_1$ denote the Hilbert space $H\times U_1$ with componentwise addition and scalar multiplication and inner product
\[
\langle (x_1,y_1),(x_2,y_2)\rangle_{H\oplus U_1}:=\langle x_1,x_2\rangle+\langle y_1,y_2\rangle_{U_1}\,\,\,\mbox{ }\forall x_1,x_2\in H\mbox{ and }\forall y_1,y_2\in U_1.
\]
The trajectory-valued random variable $(X,W):(\Omega,\mathcal{F},P)\rightarrow C([0,T];H\oplus U_1)$ is defined by
\begin{equation}\label{pe55}
(X,W)(\omega):= t\in[0,T]\mapsto(X(t)(\omega),W(t)(\omega))\,\,\,\mbox{ }\forall\omega\in\Omega
\end{equation}
and $(Y,V):(\Omega,\mathcal{F},P^\prime)\rightarrow C([0,T];H\oplus U_1)$ is defined analogously.
\begin{pro}\label{ppA1}
The trajectory-valued random variables $(X,W)$ and $(Y,V)$ have the same distribution.
\end{pro}
\textbf{Scheme of the proof.} This result is a byproduct of the proof of Theorem~\ref{ptA1}. The existence of the processes $(X(t))_{t\in[0,T]}$ and $(Y(t))_{t\in[0,T]}$ appearing in equations (\ref{e3.73aa}) and (\ref{e3.73bb}), respectively, is established using contraction maps on suitable Banach spaces of $p$-integrable processes, where $2<p<\infty$ (see \cite[Theorem 7.4]{DPZ}). We define processes $X_0(t)=x$ and $Y_0(t)=x$ for all $t\in[0,T]$. Clearly $\mathcal{L}(X_0,W)=\mathcal{L}(Y_0,V)$. For each $m\in\mathbb{N}$ let $(X_m(t))_{t\in[0,T]}$ and $(Y_m(t))_{t\in[0,T]}$ be the $m$th iterates of the respective contraction maps. One can show by induction that $\mathcal{L}(X_m,W)=\mathcal{L}(Y_m,V)$ for all non-negative integers $m$ and from this and convergence to the fixed points we get $\mathcal{L}(X,W)=\mathcal{L}(Y,V)$.

\bigskip

\noindent We conclude this appendix with a useful result which is not entirely obvious. It shows the relationship between It$\hat{\mathrm{o}}$ integrals with respect to two Wiener processes defined on related probability spaces. Let $\phi\in L^2([0,1];U)$ and let $\epsilon\in(0,1]$. Define the probability measure $P^\epsilon$ on $(\Omega,\mathcal{F})$ as in equation~(\ref{peafter24}) and define the Wiener process $(W^\epsilon(t))_{t\in[0,1]}$ on probability space $(\Omega,\mathcal{F},P^\epsilon)$ as in equation~(\ref{pe20}). In the following lemma $\mathcal{P}_1$ denotes the $(\mathcal{F}_t)$-predictable $\sigma$-algebra of subsets of $[0,1]\times\Omega$.

\begin{lem}\label{pAl1}
Let $\Phi:([0,1]\times\Omega,\mathcal{P}_1)\rightarrow (L_2(U,H),\mathcal{B}_{L_2(U,H)})$ be a measurable function such that for some positive real number $C$
\[
\int_0^1\|\Phi(s,\omega)\|^2_{L_2(U,H)}\,ds\le C\,\,\,\,\,\mbox{ for }P\mbox{ a.e. }\omega\in\Omega.
\]
Then
\[
\int_0^1\Phi(s)\,dW^\epsilon(s)=\int_0^1\Phi(s)\,dW(s)-\epsilon^{-\frac{1}{2}}\int_0^1\Phi(s)\phi(s)\,ds\,\,\,\,\,\mbox{ }P\mbox{ a.e.}.
\]
\end{lem}
\textbf{Scheme of the proof.} The result of the lemma is immediate if $\Phi$ is an elementary process. The result for a bounded function $\Phi$ is obtained by approximating $\Phi$ in $L^2([0,1]\times\Omega,\mathcal{P}_1, dt\times dP; L_2(U,H))$ with a uniformly bounded sequence of elementary processes. For general $\Phi$, we do approximation with a sequence of bounded functions obtained by truncating $\Phi$.

\paragraph{Acknowledgements}
The author is very grateful to Ben Goldys for posing this problem and for helpful discussions and would like to thank the anonymous referee whose suggestions improved this paper.

\end{document}